\documentclass[12pt,dvips]{amsart}

\usepackage{amsmath}
\usepackage{amsthm}
\usepackage{amssymb}
\usepackage{eucal}
\usepackage{times}
\usepackage{euler}
\usepackage{eufrak}

\usepackage{hyperref}

\usepackage{pst-node}
\usepackage{pst-plot}
\usepackage{pst-eps}

\usepackage[dvips]{graphicx}

\author{Piotr \'Sniady}

\address{Institute of Mathematics \\
University of Wroclaw \\ pl.~Grunwaldzki~2/4 \\ 50-384 Wroclaw
\\ Poland}

\email{Piotr.Sniady@math.uni.wroc.pl}

\title[Asymptotics of characters and genus expansion]%
{Asymptotics of characters of symmetric groups, genus expansion and
free probability}

\numberwithin{equation}{section} \numberwithin{figure}{section}

\theoremstyle{plain}
\newtheorem{lemma}{Lemma}[section]
\newtheorem{theorem}[lemma]{Theorem}
\newtheorem{proposition}[lemma]{Proposition}
\newtheorem{corollary}[lemma]{Corollary}

\newtheorem{claim}[lemma]{Claim}

\theoremstyle{definition}

\theoremstyle{remark}
\newtheorem*{remark}{Remark}

\newcommand{\A}{{\mathfrak{A}}}
\newcommand{\B}{{\mathfrak{B}}}
\newcommand{\aaa}{{\mathbf{a}}}
\newcommand{\bbb}{{\mathbf{b}}}
\newcommand{\ppp}{{\mathbf{p}}}
\newcommand{\qqq}{{\mathbf{q}}}
\newcommand{\rrr}{{\mathbf{r}}}
\newcommand{\E}{{\mathbb{E}}}
\newcommand{\C}{{\mathbb{C}}}
\newcommand{\R}{{\mathbb{R}}}
\newcommand{\Z}{{\mathbb{Z}}}
\newcommand{\M}{{\mathcal{M}}}
\newcommand{\N}{{\mathbb{N}}}
\newcommand{\gwia}{^{\star}}
\newcommand{\Sn}[1]{{{S}}_{#1}}
\newcommand{\xtra}{\ast}

\newcommand{\cytat}[1]{}

\newcommand{\fpsas}{\Phi^{\ps}_{\adm}}
\newcommand{\fasps}{\Phi^{\adm}_{\ps}}

\newcommand{\fass}{\Phi^{\adm}_{\PartPerm{A}}}
\newcommand{\fpss}{\Phi^{\ps}_{\PartPerm{A}}}
\newcommand{\fpps}{\Phi^{\pp}_{\C(\PartPerm{A})}}

\newcommand{\fppsinf}{\Phi^{\pp}_{\C(\PartPerm{\infty})}}

\newcommand{\fassSym}{\Phi^{\adm}_{\Sn{A}}}
\newcommand{\fpssSym}{\Phi^{\ps}_{\Sn{A}}}
\newcommand{\fppsSym}{\Phi^{\pp}_{\Sn{A}}}

\newcommand{\fppps}{\Phi^{\pp}_{\ps}}
\newcommand{\MJM}{M^{\JM}}
\newcommand{\Mpp}{M^{\pp}}
\newcommand{\RJM}{R^{\JM}}
\newcommand{\Rpp}{R^{\pp}}

\newcommand{\Sdwa}{{\mathbb{S}}^2}
\newcommand{\Ss}{\mathsf{S}}
\newcommand{\PartPerm}[1]{\widetilde{\Sn{#1}}}

\newcommand{\myalgebra}{\mathfrak{P}}

 \DeclareMathOperator{\comp}{comp}
 \newcommand{\compinv}{{\comp^{-1}}}
 \DeclareMathOperator{\fat}{fat}
 \DeclareMathOperator{\JM}{JM}
 \DeclareMathOperator{\adm}{as}
 \DeclareMathOperator{\ps}{ps}
 \DeclareMathOperator{\pp}{P}
 
 \DeclareMathOperator{\Tr}{Tr}

 \DeclareMathOperator{\PP}{P}
 
 \DeclareMathOperator{\NC}{NC}
 \DeclareMathOperator{\NCP}{NC_2}
 \DeclareMathOperator{\Moeb}{Moeb}

 \DeclareMathOperator{\genus}{genus}

\makeindex

\begin{document}

\begin{abstract}
The convolution of indicators of two conjugacy classes on the
symmetric group $\Sn{q}$ is usually a complicated linear combination
of indicators of many conjugacy classes. Similarly, a product of the
moments of the Jucys--Murphy element involves many conjugacy classes
with complicated coefficients. In this article we consider a
combinatorial setup which allows us to manipulate such products
easily: to each conjugacy class we associate a two-dimensional
surface and the asymptotic properties of the conjugacy class depend
only on the genus of the resulting surface. This construction
closely resembles the genus expansion from the random matrix theory.
As the main application we study irreducible representations of
symmetric groups $\Sn{q}$ for large $q$. We find the asymptotic
behavior of characters when the corresponding Young diagram rescaled
by a factor $q^{-1/2}$ converge to a prescribed shape. The character
formula (known as the Kerov polynomial) can be viewed as a power
series, the terms of which correspond to two-dimensional surfaces
with prescribed genus and we compute explicitly the first two terms,
thus we prove a conjecture of Biane.
%

\end{abstract}


\maketitle





%



\section{Introduction}

\subsection{Irreducible representations of large symmetric groups}

\label{subsec:need}

Irreducible representations of symmetric groups are in the
one-to-one correspondence with Young diagrams and due to algorithms
such as Murnaghan--Nakayama rule or Robinson--Littlewood rule the
essential questions concerning representations and characters of
symmetric groups can be answered by a combinatorial study of the
corresponding Young diagrams \cite{FultonYoungtableaux}. However,
when one studies the asymptotic properties of large symmetric
groups, the work with Young tableaux becomes cumbersome and one is
in the need to find another object which would encode the same
information in a more convenient way. Alternatively, one can state
the above problem as follows: a typical partition of a large number
$q$ (or equivalently, a Young diagram with $q$ boxes) is a
collection of at least $\sqrt{q}$ numbers, hence it contains a lot
of information. Nevertheless, we can expect that one does not need
to know all this information to extract (with a reasonable accuracy)
the properties of characters and representations. Therefore the
question arises how to compress the information about the Young
diagrams in the most efficient way.

It turns out that the right tool for development of the above
program is the notion of the transition measure of a Young diagram
which was introduced by Kerov
\cite{Kerov1993transition,Kerov1999differential}.
The transition measure of a Young diagram is a probability measure
on the real line---therefore it has a much more analytic nature than
a combinatorial notion of a Young diagram and for this reason it is
very appealing for our purposes. Yet another advantage of this
object is that it is possible to characterize it in many ways
\cite{Biane1998,Okounkov2000randompermutations}.

Of course, a question arises how to relate the transition measure of
a Young diagram with the values of the corresponding characters.
In this article we find the appropriate formula in a form of an
asymptotic series the terms of which correspond to two-dimensional
surfaces with a prescribed genus.

The first result in this direction was obtained by Biane
\cite{Biane1998}: he showed that if the sequence of Young diagrams
(after appropriate scaling) converges to some limit shape then the
leading term corresponds to surfaces with genus zero and therefore
this leading term can be computed by the means of Voiculescu's free
probability theory \cite{VoiculescuDykemaNica,HiaiPetz2000}.
Okounkov \cite{Okounkov2000randompermutations} was studying the
distribution of the length of the first rows of a large random Young
diagram sampled according to the Plancherel measure; he showed that
the limit distribution coincides with the limit distribution of the
biggest eigenvalues of a random matrix in the Gaussian Unitary
Ensemble (GUE). The main idea of his proof was the observation that
the both the formula for the moments of a GUE random matrix and the
formula for the moments of the transition measure of a random Young
diagram can be viewed as series the terms of which correspond to
two-dimensional surfaces with prescribed genus. One can view results
of this article as an attempt to simplify some of the arguments of
Biane \cite{Biane1998} and simultaneously to provide better
asymptotic expansion. Our results are closely related to those of
Okounkov \cite{Okounkov2000randompermutations} with the difference
that we do not study the connection with random matrices but on the
bright side we do not restrict ourselves to the case of the
Plancherel measure.

In the remaining part of the introduction we will present a more
detailed view of the methods and the results of this article.



\subsection{The main technical problem: convolution in $\Sn{q}$}
\label{subsec:convolutioninSq}

Convolution of central functions $f,g\in\C(\Sn{q})$ has a very
simple structure if we write them as linear combinations of
characters: if
$$f(\pi)=\sum_{\lambda\vdash q}
a_\lambda \frac{\chi^{\lambda}(e) \chi^\lambda(\pi)}{q!}, \qquad
g(\pi)=\sum_{\lambda\vdash q} b_\lambda \frac{\chi^{\lambda}(e)
\chi^\lambda(\pi)}{q!}$$ then
$$(f g)(\pi)=
\sum_{\lambda\vdash q} a_\lambda b_\lambda \frac{\chi^{\lambda}(e)
\chi^\lambda(\pi)}{ q!},$$ where the product $f g$ denotes a product
of elements of the group algebra, i.e.~the convolution of functions.
However, in many cases we cannot afford the luxury of using the
character expansion, for example in the case when we actually try to
find the asymptotics of characters. For this reason we should find
some other families of central functions on $\Sn{q}$ for which the
convolution would have a relatively simple form.

In Section \ref{subsec:definicjasigma} we shall define such a family
$\Sigma_{k_1,\dots,k_m}\in\C(\Sn{q})$ which has a particularly
simple structure: $\Sigma_{k_1,\dots,k_m}$ is (up to a normalizing
factor) an indicator of the conjugacy class of permutations with a
prescribed cycle decomposition $k_1,\dots,k_m$. This great
simplicity has very appealing consequences: it will be very easy to
evaluate functions $\Sigma$ on permutations and also it will be very
simple to write any central function on $\Sn{q}$ as a linear
combination of functions $\Sigma$ (operations which are somewhat
cumbersome for characters). This object was introduced and studied
by Ivanov and Kerov \cite{IvanovKerov1999}.

A question arises how to express a product
$\Sigma_{k_1,\dots,k_m}\cdot \Sigma_{l_1,\dots,l_n} \in \C(\Sn{q})$
as a linear combination of some other normalized conjugacy class
indicators $\Sigma$. This problem is very closely related to
calculation of the, so called, connection coefficients
\cite{GouldenJackson1996connection,GouldenJacksonLatour,
Goupil1990,Goupil1994,GoupilSchaeffer1998}. The latter problem asks
for the number of solutions of the equation $\pi_1 \pi_2 \pi_3=e$
where $\pi_1,\pi_2,\pi_3\in \Sn{q}$ must have a prescribed cycle
structure. The formulas for the connection coefficients are
available in many concrete cases, however they are not satisfactory
for the purpose of this article.

\subsubsection{The first main result: calculus of partitions and genus expansion}
It turns out that the solution to the above problem can be obtained
by considering some more general objects. In Section
\ref{subsec:pierwszawzmiankaofppsq} we shall define normalized
conjugacy class indicators $\Sigma_{\pi}$ which are indexed no
longer by sequences of integers but by partitions of finite ordered
sets. We also equip partitions with an explicit multiplicative
structure in such a way that $\Sigma$ becomes a homomorphism. For
this reason we can in fact forget in applications about the
symmetric group $\Sn{q}$ and perform all calculations in the
partition language.

In applications we are interested in the asymptotic behavior of the
contribution of various conjugacy classes when $q\to\infty$. It
turns out that the order of such a contribution of $\Sigma_{\pi}$
can be read directly from the genus of a two-dimensional surface
associated to the partition $\pi$; in this way our description is
very similar to the genus expansion for random matrices
\cite{Zvonkin1995matrixintegrals}. Even more striking is that the
degree $q$ of the symmetric group $\Sn{q}$ does not enter into the
multiplicative structure of partitions and therefore the calculus of
partitions is able to provide statements about symmetric groups
which do not depend on $q$. Another great advantage is that the
moments of the Jucys--Murphy element (which are closely related to
the moments of the transition measure) can be easily expressed
within our calculus.
All these features make the calculus of partition and its genus
expansion a perfect tool for the study of symmetric groups.

Unfortunately, one of our ultimate goals---formulas which relate
characters and the moments of the transition measure---turn out to
be quite involved and we need also some other tools to solve this
problem. As we shall see in the following, these tools are provided
by the free probability theory.


\subsection{Free probability and free cumulants}
\label{subsec:freeproba}

In this paper we study free cumulants of the transition measure of a
given Young diagram. The notion of free cumulants plays a
fundamental role in the free probability, a theory which was
initiated by Voiculescu in order to answer some old questions in the
theory of operator algebras but it soon evolved into an exciting
self--standing theory with many links to other fields, see
\cite{HiaiPetz2000,VoiculescuDykemaNica,VoiculescuPlenar,VoiculescuLectures}.
This theory can be viewed as a highly non--commutative probability
theory in which the notion of independence of random variables was
replaced by a non--commutative notion of freeness.  For the purpose
of this article we shall concentrate on the combinatorial aspect of
this theory connected with non--crossing partitions \cite{Kreweras}
and mentioned above free cumulants, which were introduced by
Speicher \cite{Speicher1994,Speicher1997,Speicher1998}.

\subsubsection{Free convolution and free cumulants}
\label{subsec:convolution} For probability measures $\mu$, $\nu$ on
the real line one can define their free convolution $\mu\boxplus\nu$
which is also a probability measure on the real line. One of the
first problems of free probability was to study this convolution.
The simplest approach is to consider the sequence of moments
\begin{equation}
M_i(\mu)=\int_{\R} x^i \ d\mu(x) \label{eq:momenty}
\end{equation}
of a given measure $\mu$ and ask what is the relation between the
moments of $\mu\boxplus\nu$ and the moments of $\mu$ and $\nu$.
Unfortunately, in turns out that the answer is given by a sequence
of quite complicated polynomials.
The solution to the above problem of finding a nice description of
the free convolution is given by free cumulants. To the measure
$\mu$ we assign a sequence of its free cumulants
$R_1(\mu),R_2(\mu),\dots$;
every free cumulant $R_i(\mu)$ is a certain polynomial in the
moments of $\mu$, and conversely, every moment of $\mu$ can be
expressed as a certain polynomial in free cumulants; therefore the
sequence of moments and the sequence of free cumulants carry the
same information. The advantage of the notion of free cumulants over
the notion of moments is the simplicity of the relation between free
cumulants: $R_n(\mu\boxplus\nu)=R_n(\mu)+R_n(\nu)$. In this article
we are not interested in the study of free convolution; our point is
that free cumulants have a miraculous property of simplifying
certain complicated non--commutative relations.

\subsubsection{Free cumulants of the transition measure}
\label{subsec:kerov}
In this article we study the relation between the transition measure
$\mu^{\lambda}$ of a Young diagram $\lambda$ and the characters of
the corresponding irreducible representation. The simplest idea
would be to describe the transition measure in terms of its moments
$M_i(\mu^{\lambda})$, however---as we already mentioned---the
relation between these moments and characters turns out to be quite
complicated. Similarly as in the example from the above Section
\ref{subsec:convolution}, free cumulants can simplify dramatically
the complexity of the formulas: it was shown by Biane
\cite{Biane1998} that for each $n$ the value $\Sigma_n$ of the
normalized character on the $n$--cycle
can be expressed as a certain polynomial in free cumulants
$R_2,R_3,\dots$ of the transition measure of the corresponding Young
diagram. Furthermore, the leading term is particularly simple,
namely
\begin{equation}
\label{eq:asymptotyka0}
\Sigma_{n}=
R_{n+1} + \text{lower degree terms}.
\end{equation}
The fundamental property of this polynomial, called \index{Kerov
polynomial} Kerov polynomial, is that it is universal for all Young
diagrams.

Kerov polynomials seem to have very interesting combinatorial
properties but not too much about them is known
\cite{BianeCharacters,StanleyRectangularCharacters,
StanleySlajdyKerov}. Kerov \index{Kerov's conjecture} conjectured
that all coefficients of Kerov polynomials are non--negative
integers; this conjecture seems to be quite difficult.
Unfortunately, Kerov's conjecture does not seem to have interesting
applications, but on the other hand its proof might be much more
interesting than the conjecture itself: Biane \index{Biane's
conjecture} \cite{BianeCharacters} suggested that the coefficients
of Kerov polynomials might be interpreted as the number of certain
intervals in the decomposition of the Cayley graph of the symmetric
group and it would be very interesting to state Biane's conjecture
in a more concrete form.

\subsubsection{The second main result: second-order expansion for Kerov polynomials}
One of the main results of this paper is a more precise asymptotic
expansion of characters
\begin{multline}
\Sigma_{n}=
  R_{n+1}+ \\  \sum_{\substack{m_2,m_3,\dots\geq 0 \\
                  2m_2+3 m_3+4 m_4+\cdots=n-1}}\!\!
            \frac{1}{4} \binom{n+1}{3} \binom{m_2+m_3+\cdots}{m_2,m_3,\dots}
                    \prod_{s\geq 2} \big( (s-1) R_s \big)^{m_s} +\\ \text{lower degree terms},
\label{eq:asymptotyka1}
 \end{multline}
which was conjectured by Biane \cite{BianeCharacters}. In other
words: we calculate explicitly the coefficients of Kerov polynomials
corresponding to the two highest--degree terms. We also outline an
algorithm which can provide such an expansion of any order.

After the first version of this article was made public
\cite{Sniady2003HigherOrder,Sniady2003pushing} Goulden and Rattan
\cite{GouldenRattan05} using different methods found an explicit
formula for all coefficients of Kerov polynomials. The Kerov's
positivity conjecture remains open until now.

\subsection{Applications: Fluctuations of random Young diagrams}

The methods presented in this article are very useful in the study
of the asymptotic properties of symmetric groups. An example of such
an application is presented in our subsequent work
\cite{Sniady2005GaussuanFluctuationsofYoungdiagrams} where we study
the distribution of a random Young diagram contributing to a given
reducible representation of the symmetric group $S_q$ in the limit
$q\to\infty$. We prove that for a large class of such
representations the fluctuations of the Young diagram around the
limit shape are asymptotically Gaussian. Our main tool in
\cite{Sniady2005GaussuanFluctuationsofYoungdiagrams} is the calculus
of partitions introduced in this article.

\subsection{Overview of the article}
In Section \ref{sec:preliminaries} we introduce the main actors: the
normalized indicators of conjugacy classes $\Sigma_{k_1,\dots,k_l}$
and the Jucys--Murphy element $J$. We also outline briefly how
important properties of Young diagrams and representations of
$\Sn{q}$ are encoded by the distribution of $J$. It will be
convenient to use in this article  the language of the probability
theory (`distribution' and `moments' of `random variables')
therefore we also introduce the necessary conventions, however a
reader might easily translate all statements into her/his favorite
language. We also recall briefly some notions connected with
partitions of finite ordered sets: fat partitions and non--crossing
partitions.


In the central part of this paper, Section \ref{sec:calculus}, we
introduce the calculus of partitions and study its properties. This
part is written as a `user--friendly user guide': we postponed all
technical and boring proofs to Section \ref{sec:proofs} in order to
allow the readers to use the calculus of partitions without
troubling why the machinery works.


%
%

In Section \ref{sec:mainresult} we study general properties of free
cumulants of the Jucys--Murphy element; in particular in Section
\ref{sec:evercrossing} we find explicitly the second--order
asymptotic expansion of these cumulants.


Section \ref{sec:proofs} is devoted to proofs of some technical
results.

In Section \ref{sec:randommatrices} we present some final remarks.
Especially interesting are Sections \ref{subsec:bianeconnection} and
\ref{subsec:okounkov} where we present connections with the work of
Biane \cite{Biane1998} and Okounkov
\cite{Okounkov2000randompermutations}. Section \ref{subsec:okounkov}
provides a natural geometric interpretation of pushing partitions
and it can be used to rephrase the results of Okounkov in a more
combinatorial language.

%
%
%
%
%
%



\section{Preliminaries}
\label{sec:preliminaries}

\cytat{
\begin{quotation}
{\it Ko\'n, jaki jest, ka\.zdy widzi.} \\ \textsc{Benedykt
Chmielowski, `Nowe Ateny albo Akademia wszelkiej sciencyi pe\l{}na'}
\end{quotation}
}

\subsection{Symmetric group}
There are many equivalent definitions of the transition measure of a
Young diagram
\cite{VershikKerov1995asymptoticbehavior,Kerov1993transition,
Kerov1999differential,OkounkovVershik1996,
Biane1998,Okounkov2000randompermutations} and for the sake of
completeness we shall recall them in the following. Nevertheless, we
shall use in this article only the description from Section
\ref{subsec:transition1}.

\subsubsection{Transition measure of a Young diagram---the
eigenvalues approach} \index{transition measure of a Young diagram}
\label{subsec:biane} The following description of the transition
measure is due to Biane \cite{Biane1998} and probably it is the
simplest one.

Consider an element $\Gamma\in \M_{q+1}\big( \C(\Sn{q}) \big)$ given
by
\begin{equation} \label{eq:macierzjucysia}
\Gamma=\left[
\begin{matrix}   0 & (1,2) & (1,3) & \dots & (1,q) & 1 \\
 (1,2) & 0 & (2,3) & \dots  & (2,q) & 1 \\
(1,3) & (2,3) & 0  & \dots  & (3,q) & 1 \\
 \vdots & \vdots & \vdots & \ddots  & \vdots & \vdots \\
 (1,q) & (2,q) & (3,q) & \dots  & 0 & 1 \\
 1 & 1 & 1 & \dots & 1 & 0
 \end{matrix} \right],
 \end{equation}
where $(i,j)\in\Sn{q}$ denotes the transposition exchanging $i$ and
$j$.

Let $\rho^{\lambda}:\C(\Sn{q})\rightarrow \M_k(\C)$ be an
irreducible representation of $\Sn{q}$ corresponding to the Young
diagram $\lambda$. We apply map $\rho^{\lambda}$ to every entry of
the matrix $\Gamma\in \M_{q+1}(\Sn{q})$ and denote the outcome by
$\rho^{\lambda}(\Gamma)\in \M_{q+1}\big( \M_{k}(\C) \big)= \M_{(q+1)
k} (\C)$. Alternatively, if we treat $\Gamma$ as an element of
$\M_{q+1}(\C) \otimes \C(\Sn{q})$ then $\rho^{\lambda}(\Gamma)= (1
\otimes \rho^{\lambda}) \Gamma\in \M_{q+1}(\C) \otimes
\M_k(\C)=\M_{(q+1)k}(\C)$.

Let $\zeta_1,\dots,\zeta_{(q+1)k}\in\R$ be the eigenvalues of the
matrix $\rho^{\lambda}(\Gamma)\in \M_{(q+1) k} (\C)$; then the
transition measure of the Young diagram $\lambda$ is the probability
measure on $\R$ which (up to a normalization) is the counting
measure of eigenvalues of $\rho^{\lambda}(\Gamma)$:
$$\mu^{\lambda}=\frac{\delta_{\zeta_1}+ \cdots + \delta_{\zeta_{(q+1)k}}}{(q+1)k}.$$

\subsubsection{Generalized Young diagrams} \index{generalized Young
diagram} \index{Young diagram!generalized}
%

\begin{figure}[tb]
\includegraphics{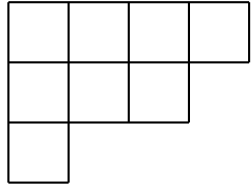}
\caption[Example of a Young diagram]{Young diagram associated with a
partition $8=4+3+1$.} \label{fig:young1}
\end{figure}

\begin{figure}[tb]
\includegraphics{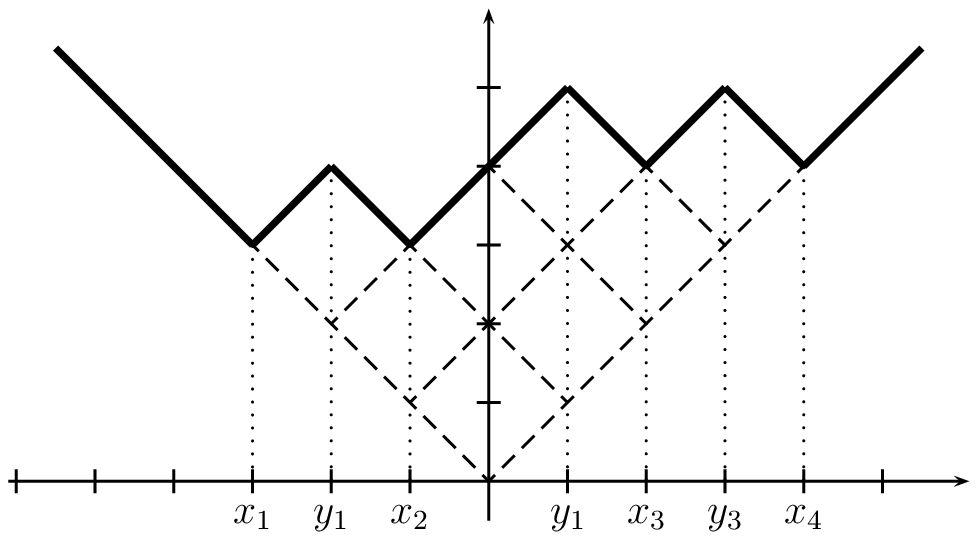}
\caption[Example of a generalized Young diagram]{Generalized Young
diagram associated with a partition $8=4+3+1$.} \label{fig:young2}
\end{figure}

Let $\lambda$ be a Young diagram. We assign to it a piecewise affine
function $\omega^\lambda:\R\rightarrow\R$ with slopes $\pm 1$, such
that $\omega^\lambda(x)=|x|$ for large $|x|$ as it can be seen on
the example from Figure \ref{fig:young2}. By comparing Figure
\ref{fig:young1} and Figure \ref{fig:young2} one can easily see that
the graph of $\omega^\lambda$ can be obtained from the graphical
representation of the Young diagram by an appropriate mirror image,
rotation and scaling by the factor $\sqrt{2}$. We call
$\omega^\lambda$ the generalized Young diagram associated with the
Young diagram $\lambda$ \cite{Kerov1993transition,
Kerov1998interlacing,Kerov1999differential}. Alternatively, we can
encode the Young diagram $\lambda$ using the sequence of local
minima of $\omega^\lambda$ (denoted by $x_1,\dots,x_m$) and the
sequence of local maxima of $\omega^\lambda$ (denoted by
$y_1,\dots,y_{m-1})$, which form two interlacing sequences of
integers \cite{Kerov1998interlacing}.

The class of generalized Young diagrams consists of all functions
$\omega:\R\rightarrow\R$ which are Lipschitz with constant $1$ and
such that $\omega(x)=|x|$ for large $|x|$ and of course not every
generalized Young diagram can be obtained by the above construction
from some Young diagram $\lambda$.

The setup of generalized Young diagrams is very useful in the study
of the asymptotic properties since it allows us to define easily
various notions of convergence of the Young diagram shapes.

\subsubsection{Transition measure---the analytic approach}
\index{transition measure of a Young diagram}
\label{subsec:transitionanalytic} To any generalized Young diagram
$\omega$ we can assign the unique probability measure $\mu^{\omega}$
on $\R$, called transition measure of $\omega$, which fulfills
\begin{equation}
\label{eq:definicja3} \log \int_{\R} \frac{1}{z-x}
d\mu^{\omega}(x)=- \frac{1}{2} \int_{\R} \log (z-x) \omega''(x) dx
\end{equation}
for every $z\notin\R$
\cite{OkounkovVershik1996,Biane1998,Okounkov2000randompermutations}.
A great advantage of this definition is that after applying
integration by parts one can easily see that the map
$\omega\mapsto\mu^{\omega}$ is continuous in many reasonable
topologies.

The generalized Young diagram $\omega^{p \lambda}:x \mapsto p
\omega^{\lambda}\big( \frac{x}{p} \big)$ corresponds to the Young
diagram $\lambda$ geometrically scaled by factor $p>0$; it is easy
to see that \eqref{eq:definicja3} implies that the corresponding
transition measure $\mu^{p \lambda}$ is a dilation of
$\mu^{\lambda}$:
\begin{equation}
\label{eq:skalowanie} \mu^{p \lambda}= D_p \mu^{\lambda}.
\end{equation}

The above definition \eqref{eq:definicja3} becomes simpler in the
case of the (usual) Young diagrams
$$ \int_{\R} \frac{1}{z-x} d\mu^{\lambda}(x) =
\frac{\prod_{1\leq i \leq m-1} (z-y_i)}{\prod_{1\leq i\leq m}
(z-x_i)} $$ and implies that the transition measure is explicitly
given by
$$\mu^\lambda = \sum_{1\leq k\leq n}
\frac{\prod_{1\leq i\leq n-1} (x_k-y_i)  }{\prod_{i\neq k}
(x_k-x_i)} \delta_{x_k}. $$


\subsubsection{Transition measure---the representation theoretic approach}
\index{transition measure of a Young diagram} It is possible to find
another interpretation of transition measure in the language of
representation theory
\cite{OkounkovVershik1996,Biane1998,Okounkov2000randompermutations}:
let us consider the representation
$[\lambda]\uparrow_{\Sn{q}}^{\Sn{q+1}}$ of $\Sn{q+1}$ induced from
the representation $[\lambda]$ of $\Sn{q}$. By the branching rule
$[\lambda]\uparrow_{\Sn{q}}^{\Sn{q+1}}$ decomposes as a direct sum
of representations corresponding to the Young diagrams obtained from
$\lambda$ by adding one box. It is possible to add a box exactly at
the minima $x_k$. The measure $\mu_\lambda$ assigns a mass to each
point $x_k$ which is proportional to the dimension of the
irreducible representation of $\Sn{q+1}$ corresponding to the
diagram $\lambda$ augmented in point $x_k$.

\subsubsection{Partial permutations}
The following notion was introduced and studied by Ivanov and Kerov
\cite{IvanovKerov1999}. \index{partial permutation}
\index{permutation!partial} A partial permutation of the set $A$ is
a pair $\alpha=(d,w)$, where $d\subseteq A$ and $w:A\rightarrow A$
is any bijection such that $d(x)=x$ for every $x\in A\setminus d$.
Set $d$ is called the support of $\alpha$. The set of partial
permutations of $A$ will be denoted by \index{$S_A$@$\PartPerm{A}$}
$\PartPerm{A}$. Given two partial permutations $(d_1,w_1)$,
$(d_2,w_2)$ we consider their product $(d_1,w_1)\cdot
(d_2,w_2):=(d_1 \cup d_2, w_1 w_2)$; given this multiplication the
set $\PartPerm{A}$ of partial permutations becomes a semigroup. By
$\Sn{A}$ we denote the permutation group of the set $A$. There is an
important homomorphism of semigroups $\PartPerm{A}\rightarrow\Sn{A}$
given by forgetting the support $(d,w)\mapsto w$; therefore every
partial permutation can be regarded as a (usual) permutation as
well.

For $A\subseteq B$ there is a homomorphism
$\theta^B_A:\C(\PartPerm{B})\rightarrow \C(\PartPerm{A})$ of partial
permutation algebras given by
$$ \theta^B_A (d,w)=\begin{cases} (d,w) & \text{if } d
\subseteq A, \\ 0 & \text{otherwise}.\end{cases}$$ By
$\C(\PartPerm{\infty})$ we denote the projective limit of algebras
$\C(\PartPerm{\{1,\dots,q\}})$ with respect to the morphisms
$\theta$.

\subsubsection{Abelian algebras in the language of probability theory}

Usually as a primary object of the probability theory one considers
a Kolmogorov space $(\Omega,{\mathcal{M}},P)$---where $\Omega$ is a
set, ${\mathcal{M}}$ is a $\sigma$--field of measurable sets and $P$
is a probability measure--- but equally well we may consider an
Abelian algebra $\A$ of all random variables on $\Omega$ with all
moments finite and the expected value $\E_\C:\A\rightarrow\C$. More
generally, if $\mathcal{M}'$ is a $\sigma$--subfield of
$\mathcal{M}$ we may consider an algebra $\B$ of
$\mathcal{M}'$--measurable random variables and a conditional
expectation $\E_\B:\A\rightarrow\B$. In this way the probability
theory becomes a theory of Abelian algebras $\A$ equipped with a
linear functional $\E_\C$ or, more generally, a theory of pairs of
Abelian algebras $\B\subset\A$ and maps $\E_\B:\A\rightarrow\B$.

By turning the picture around we may regard any Abelian algebra $\A$
equipped with a linear map $\E$ as an algebra of random variables
even if $\A$ does not arise from any Kolmogorov space (this
observation was a starting point of the non--commutative probability
theory \cite{MeyerQuantumprobability}). Now we can use the
probability theoretic language when speaking about $\A$: elements of
$\A$ can be called random variables and the numbers $\E(X^k)$ can be
called moments of the random variable $X\in\A$. Similarly we shall
interpret pairs of Abelian algebras $\B\subset\A$ equipped with a
map $\E:\A\rightarrow\B$.

In the classical probability theory any real--valued random variable
$X$ can be alternatively viewed as a selfadjoint multiplication
operator on the algebra $\A$ and according to the spectral theorem
can be written as an operator--valued integral $X=\int_\R z \
dQ(z)$, where $Q$ denotes the spectral measure of $X$. It is easy to
see that the distribution $\mu$ of the random variable $X$ and the
spectral measure of the operator $X$ are related by equality
$\mu(F)=\E_\C[ Q(F)]$ for any Borel set $F\subseteq\R$. This simple
observation can be used in our new setup to define a distribution
$\mu$ of an element $X=X\gwia\in\A$ to be a $\B$--valued measure on
$\R$ such that $\mu(F)=\E_\B[ Q(F)]$ for any Borel set
$F\subseteq\R$. If $X$ is bounded then its distribution can be
alternatively described by the moment formula
$$\int_\R x^n \ d\mu(x)=\E_\B(X^n).$$

\subsubsection{Symmetric group in the language of probability theory}
\label{subsec:symmetricgroupprobability} Let $A$ be a finite set and
let $\xtra$ be an extra distinguished element such that $\xtra\notin
A$. In the following, if not stated otherwise, we set
$A=\{1,2,\dots,q\}$.  Group $\Sn{A}$ will be regarded as a subgroup
of $\Sn{A\cup \{\xtra\}}$. In our case the role of the algebra $\A$
of random variables will be played by some Abelian subalgebra of
$\C(\Sn{A\cup\{\xtra\} })$ and the role of the smaller algebra $\B$
will be played by the center of $\C(\Sn{A})$. The conditional
expectation $\E_{\C(\Sn{A})} :\C(\Sn{A\cup\{\xtra\}})\rightarrow
\C(\Sn{A}) $ will be given by the orthogonal projection, i.e.
$$ \E_{\C(\Sn{A})}(\sigma)=\begin{cases} \sigma & \text{if } \sigma\in\Sn{A} \\
0 & \text{if } \sigma\notin\Sn{A} \end{cases}=
\begin{cases} \sigma & \text{if } \sigma(\xtra)=\xtra, \\
0 & \text{if } \sigma(\xtra)\neq \xtra .\end{cases}
$$
In the rest of this article, when it does not lead to confusions,
instead of $\E_{\C(\Sn{A})}$ we shall simply write
\index{$E$@$\E:\C(\Sn{A\cup\{\xtra\}})\rightarrow \C(\Sn{A}) $}
$\E$.

Suppose that some finite--dimensional representation
$\rho:\Sn{A}\rightarrow \M_k(\C)$ is given. This allows us to define
a scalar--valued expectation
$\E_\C:\C(\Sn{A\cup\{\xtra\}})\rightarrow\C$ by
\begin{equation}
\label{eq:wartoscoczekiwanahihihi} \E_\C(X)=\frac{1}{k} \Tr
\rho\big(\E_{\C(\Sn{A})}(X)\big)
\end{equation}
 and
to consider the distribution $\mu^{\rho}$ of a random variable
$X\in\C(\Sn{A\cup\{\xtra\}})$ with respect to this new expectation;
in this new setup $\mu^{\rho}$ is a usual (scalar--valued)
probability measure on $\R$.

%

\subsubsection{Jucys--Murphy element. Transition measure---the algebraic approach.
The algebra $\myalgebra$} \index{transition measure of a Young
diagram} \label{subsec:transition1} In the group algebra
$\C(\Sn{A\cup\{\xtra\}})$ we consider the Jucys--Murphy element
\index{Jucys--Murphy element} \index{$J$}
$$ J= \sum_{a\in A} (a \xtra), $$
where $(i j)$ denotes the transposition exchanging $i$ and $j$. We
define moments of the Jucys--Murphy element by
\index{$MJM$@$\MJM_k$}
\begin{equation}
\MJM_k=\E(J^k)= \sum_{a_1,\dots,a_k\in A } \E[ (a_1 \xtra) \cdots
(a_k \xtra) ] \in \C(\Sn{A}). \label{eq:wielkasuma}
\end{equation}
In Section \ref{subsec:momentmap} of we shall consider an extension
of this concept.

For a given Young diagram $\lambda\vdash |A|$ we consider
$\rho=\rho^\lambda$ to be an irreducible representation of $\Sn{A}$
corresponding to $\lambda$ and consider the distribution
$\mu^{\lambda}:=\mu^{\rho^{\lambda}}$ of the element $J$ with
respect to the expected value \eqref{eq:wartoscoczekiwanahihihi}.
We call $\mu^{\lambda}$ the transition measure of $\lambda$
\cite{VershikKerov1995asymptoticbehavior,Kerov1993transition,
OkounkovVershik1996,Biane1998,Okounkov2000randompermutations}. This
definition is very algebraic and for this reason it will be our
favorite definition of the transition measure.


We can treat the moments $\MJM_k$ of the Jucys--Murphy element as
elements of the partial permutations algebra $\C(\PartPerm{A})$; in
order to do this we treat every non--zero summand on the right--hand
side of \eqref{eq:wielkasuma} as a partial permutation with the
support $\{a_1,\dots,a_k\}$. It is easy to check that for
$A\subseteq B$ the morphism $\theta^B_A$ maps
$\MJM_k\in\C(\PartPerm{B})$ to $\MJM_k\in\C(\PartPerm{A})$ hence the
projective limit of the elements $\MJM_k\in\C(\PartPerm{q})$ exists
and will be denoted by the same symbol
$\MJM_k\in\C(\PartPerm{\infty})$.


We denote by $\myalgebra$ the algebra generated by elements
$\MJM_k\in\C(\PartPerm{\infty})$. It is easy to check that
$\myalgebra$ is commutative. \index{$P$@$\myalgebra$} Elements of
$\myalgebra$ can be also regarded as elements of $\C(\PartPerm{q})$
and $\C(\Sn{q})$.

\subsubsection{Gradation on $\myalgebra$} \label{subsec:degree1}
Since we would like to study asymptotic properties of Young
diagrams, we need to specify what kind of scaling we are interested
in.

Let a sequence of Young diagrams $(\lambda_N)$ be given,
$\lambda_N\vdash N$. It is natural to consider generalized Young
diagrams which correspond to geometrically scaled diagrams
$\frac{1}{\sqrt{N}} \lambda_N$ (observe that all such scaled
diagrams have the same area equal to $2$). Suppose that the shape of
the scaled diagrams converges in some sense toward a generalized
Young diagram $\lambda$. In many natural topologies this implies the
convergence of scaled moments of the transition measure, cf
\eqref{eq:skalowanie}
\begin{multline*}
\lim_{N\rightarrow\infty} \frac{1}{N^{k/2}} \E_\C (\MJM_k)=
 \lim_{N\rightarrow\infty} \E_{\C} \left[ \left( \frac{1}{N^{1/2}} J \right)^k \right]=
\\
\lim_{N\rightarrow\infty} \frac{1}{N^{k/2}} \int_\R x^k
d\mu^{\lambda_N}(x) =  \int_\R x^k d\mu^{\lambda}(x).
\end{multline*}
In the study of asymptotic properties of some elements of the group
algebra $\C(\Sn{N})$ we would like to group summands which have
asymptotically the same growth for large $N$ and for this reason
$\MJM_k$ can be treated as a monomial in ${\sqrt{N}}$ of degree $k$.

More formally, we consider a gradation \index{gradation} on
$\myalgebra$ by setting
\begin{equation} \label{eq:gradation} \deg \MJM_k=k.
\end{equation}
In Corollary \ref{cor:algebraicallyfree} we will show that
$\MJM_2,\MJM_3,\dots\in\C(\PartPerm{\infty})$ are algebraically free
and therefore this gradation is well defined. This gradation
coincides with the weight gradation considered by Ivanov and
Olshanski \cite{IvanovOlshanski2002}.

\subsubsection{Normalized conjugacy class indicators}
\label{subsec:definicjasigma}

Let integer numbers $k_1,\dots,k_m\geq 1$ be given. We define the
normalized conjugacy class indicator \index{normalized conjugacy
class indicator} \index{$Sigma_{k_1}$@$\Sigma_{k_1,\dots,k_m}$}
$\Sigma_{k_1,\dots,k_m}\in\C(\Sn{A})$ as follows
\cite{KerovOlshanski1994,BianeCharacters}:
\begin{equation}
\label{eq:definicjasigma} \Sigma_{k_1,\dots,k_m}=
\sum_{a} (a_{1,1},a_{1,2},\dots,a_{1,{k_1}}) \cdots
(a_{m,1},a_{m,2},\dots,a_{m,k_m}), \end{equation}
where the sum runs over all one--to--one functions
$$a:\big\{ \{r,s\}: 1\leq r\leq m, 1\leq s\leq {k_r}\big\}\rightarrow A$$ and
$(a_{1,1},a_{1,2},\dots,a_{1,{k_1}}) \cdots
(a_{m,1},a_{m,2},\dots,a_{m,k_m})$ is a product of disjoint cycles.
Of course, if $|A|<k_1+\cdots+k_m$ then the above sum runs over an
empty set and $\Sigma_{k_1,\dots,k_m}=0$. If any of the numbers
$k_1,\dots,k_m$ is equal to $0$ we set $\Sigma_{k_1,\dots,k_m}=0$.

In other words, let $k_1'\geq \cdots \geq k_m'$ be the sequence
$k_1,\dots,k_m$ sorted decreasingly; we consider a Young diagram
$(k_1',\dots,k_m')$ and all ways of filling it with the elements of
the set $A$ in such a way that no element appears more than once.
Each such a filling can be interpreted as a permutation when we
treat rows of the Young tableau as disjoint cycles.


We can also treat $\Sigma_{k_1,\dots,k_m}$ as a an element of the
partial permutations algebra $\C(\PartPerm{A})$; in order to do
this, we treat every summand $(a_{1,1},a_{1,2},\dots,a_{1,{k_1}})
\cdots (a_{m,1},a_{m,2},\dots,a_{m,k_m})$ as a partial permutation
with support $\{ a_{ij} : 1\leq i\leq m, 1\leq j\leq {k_i}  \}$. It
is easy to check that for $A\subseteq B$ the morphism $\theta^B_A$
maps $\Sigma_{k_1,\dots,k_m}\in\C(\PartPerm{B})$ to
$\Sigma_{k_1,\dots,k_m}\in\C(\PartPerm{A})$ hence the projective
limit of the elements $\Sigma_{k_1,\dots,k_m}\in\C(\PartPerm{q})$
exists and will be denoted by the same symbol
$\Sigma_{k_1,\dots,k_m}\in\C(\PartPerm{\infty})$. It is easy to
check that elements $\Sigma_{k_1,\dots,k_m}\in\C(\PartPerm{\infty})$
commute.

%
%

\subsubsection{Filtration}

We consider a filtration \index{filtration} on the commutative
algebra generated by
$\Sigma_{k_1,\dots,k_m}\in\C(\PartPerm{\infty})$ by setting
\begin{equation}
\label{eq:filtration} \deg
\Sigma_{k_1,\dots,k_m}=(k_1+1)+(k_2+1)+\cdots+(k_m+1).
\end{equation}
It is easy to check that the elements from the family
$\Sigma_{k_1,\dots,k_m}\in\C(\PartPerm{\infty})$ indexed by integer
$m\geq 0$ and $k_1\geq \cdots \geq k_m\geq 1$ are linearly
independent and therefore this definition makes sense but it is not
clear that this formula indeed defines a filtration; we shall prove
it in Corollary \ref{cor:filtration}. The reason for studying this
filtration is that---as we shall see in Theorem
\ref{theo:kontrolastopnia}---elements
$\MJM_k\in\C(\PartPerm{\infty})$ and
$\Sigma_{k_1,\dots,k_m}\in\C(\PartPerm{\infty})$ generate the same
commutative algebra denoted by $\myalgebra$ and filtration
\eqref{eq:filtration} is induced by the gradation
\eqref{eq:gradation}.

\subsection{Partitions} \label{sec:combinatorial}

\subsubsection{Partitions} \label{subsec:partitions}
\index{partition} We recall that $\pi=\{\pi_1,\dots,\pi_r\}$ is a
partition of a finite ordered set $X$ if sets $\pi_1,\dots,\pi_r$
are nonempty and disjoint and if $\pi_1\cup\cdots\cup\pi_r=X$. We
denote the set of all partitions of a set $X$ by $\PP(X)$. We call
sets $\pi_1,\dots,\pi_r$ blocks of the partition $\pi$. We say that
elements $a,b\in X$ are connected by the partition $\pi$ if they are
the elements of the same block.
We will call elements of the set $X$ the labels of the vertices of
the partition $\pi$. As usually, the numbers $|\pi_1|,\dots,|\pi_r|$
denote the numbers of elements in consecutive blocks of $\pi$.

We say that a block $\pi_s$ is trivial if it contains only one
element and a trivial partitions has all blocks
trivial\index{trivial partition}\index{partition!trivial}. We say
that a partition is a pair partition\index{pair
partition}\index{partition!pair} if all its blocks contain exactly
two elements.

We say that a partition $\pi$ is finer than a partition $\rho$ of
the same set if every block of $\pi$ is a subset of some block of
$\rho$ and we denote it by \index{$Pi$@$\pi\leq\rho$} $\pi\leq\rho$.
The restriction of a partition $\pi$ \index{partition!restriction
of}to a set $Y\subseteq X$ is a partition of $Y$ which connects
elements $a,b\in Y$ if and only if $a,b$ are connected by $\pi$.

In the following we present some constructions on partitions of the
set $X=\{1,2,\dots,n\}$. However, it should be understood that by a
change of labels these constructions can be performed for any finite
ordered set $X$.

\begin{figure}[tb]
\includegraphics{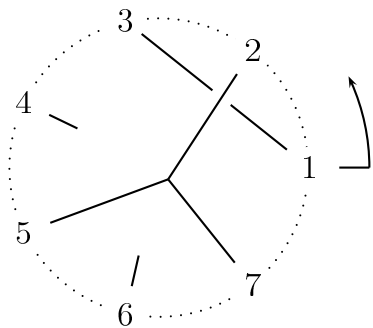}
\caption[Graphical representation of a partition]{Graphical
representation of a partition $\big\{ \{1,3\},\{2,5,7\},\{4\},\{6\}
\big\}$.} \label{fig:przecinajacapartycja}
\end{figure}

It is very useful to represent partitions graphically by arranging
the elements of the set $X$ counterclockwise on a circle and joining
elements of the same block by a line, as it can be seen on Figure
\ref{fig:przecinajacapartycja}. Note that we do not need to write
labels on the vertices in order to recover the partition from its
graphical representation if we mark the first element as the
`starting point'.

\subsubsection{Fat partitions} \index{fat partition}
\index{partition!fat} \index{$Pi$@$\pi_{\fat}$} \label{subsec:fat}
Let $\pi=\{\pi_1,\dots,\pi_r\}$ be a partition of the set
$\{1,\dots,n\}$. For every $1\leq s\leq r$ let
$\pi_s=\{\pi_{s,1},\dots,\pi_{s,l_s}\}$ with
$\pi_{s,1}<\cdots<\pi_{s,l_s}$. We define $\pi_{\fat}$, called fat
partition of $\pi$, to be a pair partition of the $2n$--element
ordered set $\{1,1',2,2',\dots,n,n'\}$ given by
$$\pi_{\fat}=
\big\{  \{ \pi_{s,t}',\pi_{s,t+1}  \} :  1\leq s \leq r \text{ and }
1 \leq t \leq l_s  \big\}, $$ where it should be understood that
$\pi_{s,l_s+1}=\pi_{s,1}$.

\begin{figure}[tb]
\includegraphics{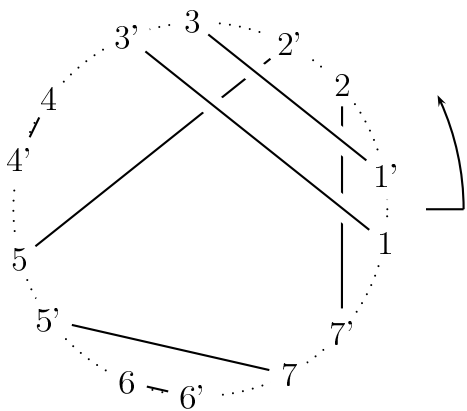}
\caption[Fat partition $\pi_{\fat}$]{Fat partition $\pi_{\fat}$
corresponding to the partition $\pi$ from Figure\
\ref{fig:przecinajacapartycja}.} \label{fig:tlustapartycja}
\end{figure}

This operation can be easily described graphically as follows: we
draw blocks of a partition with a fat pen and take the boundary of
each block, as it can be seen on Figure \ref{fig:tlustapartycja}.
This boundary is a collection of lines hence it is a pair partition.
However, every vertex $k\in\{1,\dots,n\}$ of the original partition
$\pi$ has to be replaced by its `right' and `left' copy (denoted
respectively by $k$ and $k'$). Please note that in the graphical
representation of $\pi_{\fat}$ we mark the space between $1$ and
$1'$ as the `starting point'.


\subsubsection{Non--crossing partitions} 
\begin{figure}[tb]
\includegraphics{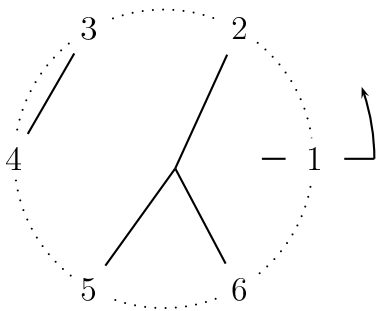}
\caption[Graphical representation of a non--crossing partition]{%
Graphical representation of a non--crossing partition $\big\{
\{1\},\{2,5,6\}, \{3,4\} \big\}$.}
\label{fig:nieprzecinajacapartycja}
\end{figure}

\index{non--crossing
partition} \index{partition!non--crossing} \index{$NC$@$\NC(X)$} We
say that a partition $\pi$ is non--crossing
\cite{Kreweras,SimionUllman,Speicher1994,Speicher1997,Speicher1998}
if for every $i\neq j$ and $a,c\in\pi_i$ and $b,d\in\pi_j$ it cannot
happen that $a<b<c<d$. For example, the partition from Figure
\ref{fig:przecinajacapartycja} is crossing (as it can be easily seen
from its graphical representation) while a partition from Figure
\ref{fig:nieprzecinajacapartycja} is non--crossing. The set of all
non--crossing partitions of a set $X$ will be denoted by $\NC(X)$
and the set of all non--crossing pair partitions of a set $X$ will
be denoted by $\NCP(X)$.

\subsubsection{Kreweras complementation map} The function
$\pi\mapsto\pi_{\fat}$ which maps partitions of a $n$--element set
to pair partitions of a $2n$--element set is one to one but in
general is not a bijection. However, it is a bijection between
$\NC(1,2,\dots,n)$ and $\NCP(1,1',\dots,n,n')$ \cite{Kreweras}.

\begin{figure}[tb]
\includegraphics{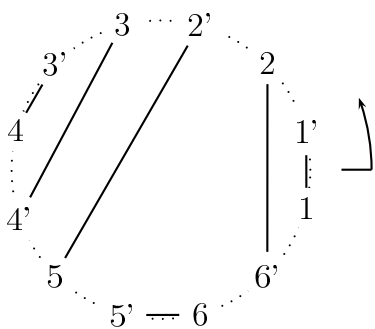}
\caption[Fat partition $\pi_{\fat}$ for a non--crossing
partition]{Fat partition $\pi_{\fat}$ for a non--crossing partition
$\pi$ from Figure \ref{fig:nieprzecinajacapartycja}.}
\label{fig:tlustapartycjanieprzecinajaca}
\end{figure}

\begin{figure}[tb]
\includegraphics{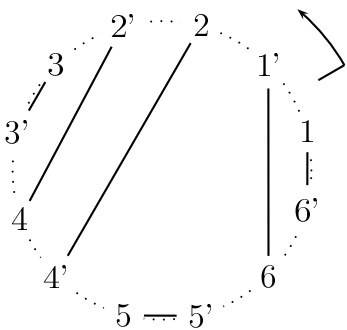}
\caption[Graphical representation of $r(\rho)$]{Graphical
representation of $r(\rho)$ where $\rho$ is given by Figure
\ref{fig:tlustapartycjanieprzecinajaca}.} \label{fig:rotacja}
\end{figure}

If $\rho$ is a pair partition of the set $\{1,1',2,2',\dots,n,n'\}$,
we denote by $r(\rho)$ a pair partition of the same set
$\{1,1',\dots,n,n'\}$ which is obtained by cyclic change of labels
$\cdots\rightarrow 2'\rightarrow 2\rightarrow 1'\rightarrow
1\rightarrow n'\rightarrow n \rightarrow\cdots$ carried by the
vertices. In the graphical representation of a partition this
corresponds to a shift of the starting point by one counterclockwise
(cf Figure \ref{fig:rotacja}). Of course $r$ is a bijection of the
set $\NCP(1,1',\dots,n,n')$.

It follows that the map $\pi\mapsto\pi_{\comp}$, called
\index{Kreweras complement} \index{$Pi$@$\pi_{\comp}$} Kreweras
complementation map, given by $(\pi_{\comp})_{\fat}=r(\pi_{\fat})$
is well--defined and is a permutation of $\NC(1,2,\dots,n)$. Due to
the canonical bijection between $\{1,2,\dots,n\}$ and
$\{1',2',\dots,n'\}$ we can identify $\NC(1,2,\dots,n)$ with
$\NC(1',2',\dots,n')$ and in the future we shall sometimes regard
$\pi_{\comp}$ as an element of $\NC(1,2,\dots,n)$ and sometimes as
an element of $\NC(1',2',\dots,n')$. It will be clear from the
context which of the options we choose.

\begin{figure}[tb]
\includegraphics{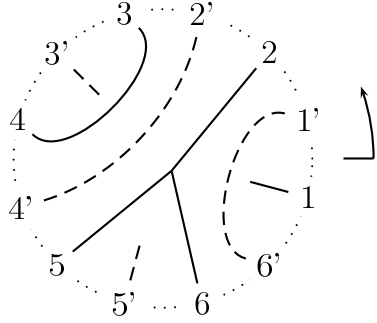}
\caption[Non--crossing partition and its Kreweras
complement]{Partition $\pi$ from Figure
\ref{fig:nieprzecinajacapartycja} in a solid line and its Kreweras
complement in a dashed line.} \label{fig:Kreweras}
\end{figure}

One can also state the above definition as follows: $\pi_{\comp}$ is
the biggest non--crossing partition of the set $\{1',2',\dots,n'\}$
with a property that $\pi\cup\pi_{\comp}$ is a non--crossing
partition of the set $\{1,1',2,2',\dots,n,n'\}$ (cf Figure
\ref{fig:Kreweras}).

We shall also consider the inverse of the Kreweras complementation
map: for $\pi, \rho\in\NC(1,2,\dots,n)$ we have
$\pi_{\compinv}=\rho$ if and only if $\pi=\rho_{\comp}$.
\index{$Pi$@$\pi_{\compinv}$} One can observe that $\rho_{\compinv}$
and $\rho_{\comp}$ are always cyclic rotations of each other.

\section{The first main result: Calculus of partitions}
\label{sec:calculus}

\cytat{
\begin{quotation}
{\it `and what is the use of a book,'} thought Alice {\it `without
pictures or conversation?'}
\\
\textsc{Lewis Carroll, `Alice's Adventures in Wonderland'}
\end{quotation}
}

\subsection[Partition--indexed conjugacy class indicator $\Sigma_\pi$]%
{Partition--indexed conjugacy class indicator $\Sigma_\pi$}

\subsubsection{Definition of\/ $\Sigma_\pi$}
\label{subsec:pierwszawzmiankaofppsq} Let $\pi$ be a partition of
the set $\{1,\dots,n\}$. Since the fat partition $\pi_{\fat}$
connects every element of the set $\{1',2',\dots,n'\}$ with exactly
one element of the set $\{1,2,\dots,n\}$, we can view $\pi_{\fat}$
as a bijection
$\pi_{\fat}:\{1',2',\dots,n'\}\rightarrow\{1,2,\dots,n\}$. We also
consider a bijection
$c:\{1,2,\dots,n\}\rightarrow\{1',2',\dots,n'\}$ given by
$\dots,3\mapsto 2', 2\mapsto 1', 1\mapsto n', n\mapsto
(n-1)',\dots$. Finally, we consider a permutation $\pi_{\fat}\circ
c$ of the set $\{1,2,\dots,n\}$.

\begin{figure}[bt]
\includegraphics{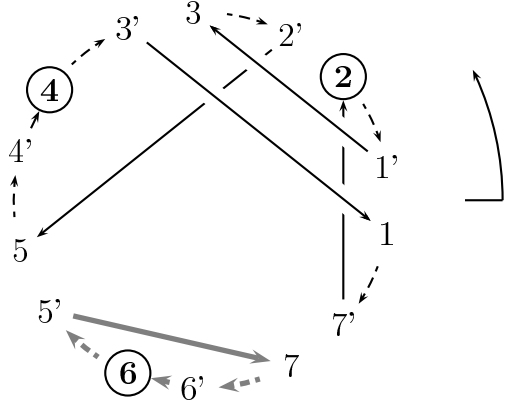}
\caption[Bijection corresponding to $\pi_{\fat}$ from Figure
\ref{fig:tlustapartycja}]{Bijection corresponding to the partition
$\pi_{\fat}$ from Figure \ref{fig:tlustapartycja} plotted with a
solid line and the bijection $c$ plotted with a dashed line. Lines
corresponding to two different cycles were plotted in different
colors (gray and black). The meaning of the additional decoration of
some of the vertices will be explained in Section
\ref{subsec:explicit}.} \label{fig:tlustapartycja2}
\end{figure}

For example, for the partition $\pi$ given by Figure
\ref{fig:przecinajacapartycja} the composition $\pi_{\fat}\circ c$
has a cycle decomposition $(1,{2},3,5,{4})({6},7)$, as it can be
seen from Figure \ref{fig:tlustapartycja2}.

Alternatively, one can identify the set $\{1',2',\dots,n'\}$ with
the set $\{1,2,\dots,n\}$; then $\pi_{\fat}$ becomes a permutation
which coincides with the construction from
\cite{Kreweras,Biane1997crossings,Biane1998} and $c$ is equal to the
full cycle $(n,n-1,\dots,2,1)$.

\begin{figure}[bt]
\includegraphics{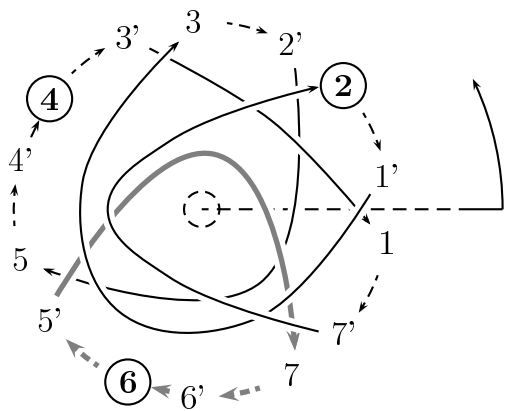}
\caption[Version of Figure \ref{fig:tlustapartycja2} in which all
lines wind clockwise]{A version of Figure \ref{fig:tlustapartycja2}
in which all lines wind clockwise around the central disc.}
\label{fig:tlustapartycja3}
\end{figure}

We decompose the permutation
$$\pi_{\fat}\circ c=(b_{1,1},b_{1,2},\dots,b_{1,j_1})
\cdots (b_{t,1},\dots,b_{t,j_t})$$ as a product of disjoint cycles.
Every cycle $b_s=(b_{s,1},\dots,b_{s,j_s})$ can be viewed as a
closed clockwise path on a circle and therefore one can compute how
many times it winds around the circle.  It might be useful to draw a
line between the central disc and the starting point (cf Figure
\ref{fig:tlustapartycja3}); the number of winds is equal to the
number of times a given cycle crosses this line and therefore is
equal to the number of indices $1\leq i\leq j_s$ such that
$b_{s,i}\leq b_{s,i+1}$, where we use the convention that
$b_{s,j_s+1}=b_{s,1}$.

To a cycle $b_s$ we assign the number
\begin{multline}
\label{eq:formulanak} k_s=(\text{number of elements in a cycle }
b_s)-\\ (\text{number of clockwise winds of }b_s).
\end{multline}

In the above example we have $b_1=(1,{2},3,5,{4})$, $b_2=({6},7)$
and $k_1=2$, $k_2=1$, as it can be seen from Figure
\ref{fig:tlustapartycja3}, where all lines clockwise wind around the
central disc.

\begin{figure}[bt]
\includegraphics{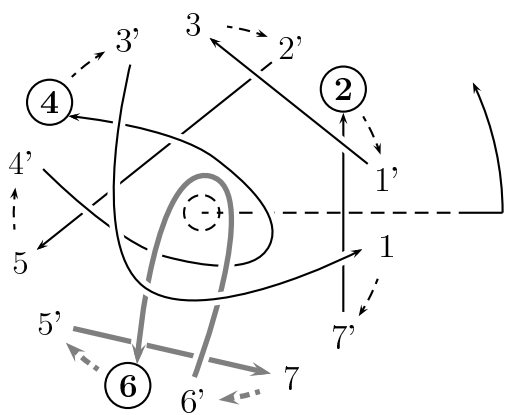}
\caption[Version of Figure \ref{fig:tlustapartycja2} in which all
lines of $\pi_{\fat}$ wind counterclockwise]{A version of Figure
\ref{fig:tlustapartycja2} in which all lines of $\pi_{\fat}$ wind
counterclockwise around the central disc.}
\label{fig:tlustapartycja4}
\end{figure}

Alternatively, we can treat every cycle $b_s$ of $\pi_{\fat} \circ
c$ as a closed counter--clockwise path on a circle (cf Figure
\ref{fig:tlustapartycja4}). Equation \eqref{eq:formulanak} implies
that
\begin{equation}
\label{eq:formulanakk} k_s=(\text{number of counterclockwise winds
of the cycle } b_s).
\end{equation}


Let a positive integer $q$ be given. In Section \ref{sec:proofs} we
shall define the most important tool of this article: a certain map
$\fppsinf$ from the set of partitions to $\C(\PartPerm{\infty})$.
The following Claim states that $\fppsinf(\pi)$ is closely related
with the normalized conjugacy class indicators from Section
\ref{subsec:definicjasigma} of Section \ref{sec:preliminaries} and
hence $\Sigma_{\pi}:=\fppsinf(\pi)$
\index{$Sigma_{\pi}$@$\Sigma_{\pi}$} can be viewed as a normalized
conjugacy class indicator which is indexed by a partition $\pi$
instead of a tuple of integers $(k_1,\dots,k_t)$.

\begin{claim}
\label{claim:definicjasigma} Let $\pi$ be a partition of the set
$\{1,2,\dots,n\}$ and let numbers $k_1,\dots,k_t$ be given by the
above construction. Then
$$\Sigma_\pi:=\fppsinf(\pi)=\Sigma_{k_1,\dots,k_t},$$
where $\Sigma_{k_1,\dots,k_t}\in\C(\PartPerm{\infty})$ on the
right--hand side should be understood as in Section
\ref{subsec:definicjasigma}.
\end{claim}
This Claim will follow from Theorem
\ref{theo:stwierdzeniepierwszejestprawda}. Reader not interested in
too technical details may even take this Claim as a very simple
alternative definition of $\Sigma_\pi=\fppsinf(\pi)$.


\subsubsection{Basic properties of $\Sigma_{\pi}$} We say that
partitions $\pi$ and $\pi'$ are cyclic rotations of each other if
one can obtain $\pi'$ from $\pi$ by a cyclic shift of labels of
vertices. Graphically, this corresponds to a change of a marked
starting point on a circle.
\begin{proposition}
\label{prop:rotacjajestok} Let $\pi$ and $\pi'$ be two partitions of
the same set which are cyclic rotations of each other. Then
$$\Sigma_{\pi}=\Sigma_{\pi'}.$$
\end{proposition}
\begin{proof}
It is enough to observe that the winding numbers considered in the
algorithm from Section \ref{subsec:pierwszawzmiankaofppsq} do not
depend on the choice of the starting point.
\end{proof}


\begin{proposition}
If partition $\pi$ connects some neighbor elements then
$\Sigma_\pi=0$.
\end{proposition}
\begin{proof}
If partition $\pi$ connects some $r$ with $r+1$ then $\pi_{\fat}$
connects $r'$ with $r+1$, hence $\pi_{\fat}(r')=r+1$ and
$(\pi_{\fat}\circ c)(r+1)=r+1$. One of the cycles of
$\pi_{\fat}\circ c$ is a $1$--cycle $b_s=(r+1)$ and one can easily
check that $k_s=0$. Claim \ref{claim:definicjasigma} implies that
$\Sigma_\pi=\Sigma_{k_1,\dots,k_t}=0$.
\end{proof}

\subsubsection{Degree of\/ $\Sigma_\pi$ and its geometric
interpretation} We keep the notation from the previous Section; in
particular $t$ denotes the number of cycles of $\pi_{\fat}\circ c$
and $r$ the number of blocks of $\pi$. Observe that
$$\sum_{1\leq s\leq t} (\text{number of elements in a cycle } b_s)=n.$$

Furthermore $\sum_{1\leq s\leq t} (\text{number of clockwise winds
of } b_s)$ is equal to the total number of clockwise winds of
$\pi_{\fat}\circ c$. When we identify sets $\{1,2,\dots,n\}$ and
$\{1',2',\dots,n'\}$ then $c$ is a full clockwise cycle and it
contributes with a one wind. It is easy to observe that every block
$\pi_s$ of partition $\pi$ is a counterclockwise cycle of the
permutation $\pi_{\fat}$, hence when we treat it as a collection of
clockwise steps it contributes with $|\pi_s|-1$ clockwise winds.  It
follows that
\begin{multline*}\sum_{1\leq s\leq t} (\text{number of clockwise winds of } b_s)=\\
1+\sum_{1\leq s\leq r} (|\pi_s|-1)=n-r+1
\end{multline*}
since $\sum_{1\leq s\leq r} |\pi_s|=n$. Therefore
\eqref{eq:formulanak} gives us
\begin{equation}
\label{eq:formulanastopien} \deg \Sigma_{\pi}=
k_1+\dots+k_t+t=r+t-1,
\end{equation}
where the degree is taken with respect to the filtration
\eqref{eq:filtration}.

\begin{figure}[tb]
\includegraphics{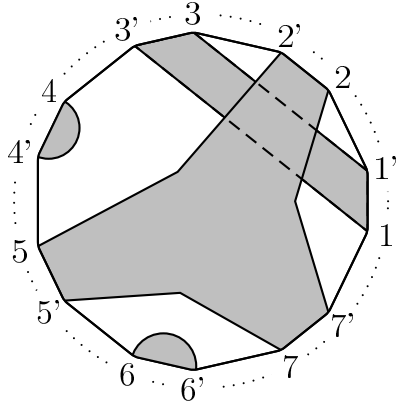}
\caption[The first collection of discs for $\pi$ from Figure
\ref{fig:przecinajacapartycja}.]{The first collection of discs for
partition $\pi$ from Figure \ref{fig:przecinajacapartycja}.}
\label{fig:firstcollection}
\end{figure}

\label{subsec:zaklejanie} We can find a natural geometric
interpretation of \eqref{eq:formulanastopien}: consider a large
sphere with a small circular hole. The boundary of this hole is the
circle that we consider in the graphical representations of
partitions. Some pairs of points on this circle are connected by
lines: the blocks of the partition $\pi_{\fat}$. To the arcs on the
boundary of the circle between points $1,1',2,2',\dots,n,n'$ and to
the lines of $\pi_{\fat}$ we shall glue two
 collections of discs. Every disc from the first collection
corresponds to one of the blocks of $\pi$ (cf Figure
\ref{fig:firstcollection}).

After gluing the first collection of discs, our sphere becomes a
surface with a number of holes. The boundary of each hole
corresponds to one of the cycles of $\pi_{\fat}\circ c$ and we shall
glue this hole with a disc from the second collection.


\index{genus of a partition} \index{partition!genus of}
\index{$Genus$@$\genus_\pi$} Thus we obtained an orientable surface
without a boundary. We define the genus of the partition $\pi$ to be
the genus of this surface. The numbers of vertices, edges and faces
of the polyhedron constructed above are the following: $V=2n$,
$E=3n$, $F=1+r+t$ therefore the genus of this surface is equal to
\begin{equation}
\label{eq:genus} \genus_\pi=\frac{2-V+E-F}{2}=\frac{n+1-r-t}{2}.
\end{equation}
We proved therefore the following result.
\begin{proposition}
\label{prop:genusowatosc} For any partition $\pi$ of an $n$--element
set
\begin{equation}
\label{eq:genusowatosc} \deg \Sigma_{\pi}=n-2 \genus_\pi,
\end{equation}
where the degree is taken with respect to the filtration
\eqref{eq:filtration}.
\end{proposition}


We leave the proof of the following simple result to the reader.
\begin{proposition}
\label{prop:genuszero} A partition is non--crossing if and only if
its genus is equal to zero.
\end{proposition}

\subsection{Multiplication of partitions}
\label{sec:mnozeniepartycji} \subsubsection{Definition of
multiplication} We will introduce an algebraic structure on the
partitions. For any finite ordered set $X$ one can consider a linear
space spanned by all partitions of $X$. We will define a
`multiplicative' structure as follows: let
$\rho=\{\rho_1,\dots,\rho_r\}$ be a non--crossing partition of a
finite ordered set $X$ and for every $1\leq s\leq r$ let $\pi^s$ be
a partition of the set $\rho_s$. We define the `$\rho$--ordered
product' of partitions $\pi^s$ by \index{$Rho$@$\rho$--ordered
product} \index{$Pi_s$@$\prod_{s} \pi^s$}
\begin{equation}
\label{eq:definicjamnozenia} \prod_{s} \pi^s :=\sum_{\sigma} \sigma,
\end{equation}
where the sum denotes a formal linear combination and it runs over
all partitions $\sigma$ of the set $X$ such that
\begin{enumerate}
\item for any $a,b\in\rho_s$, $1\leq s\leq r$ we have that $a$ and
$b$ are connected by $\sigma$ if and only if they are connected by
$\pi^s$, \item $\sigma\geq \rho_{\compinv}$,
\end{enumerate}
where $\geq$ denotes the order on partitions we introduced in
Section \ref{subsec:partitions}.

Above we defined the product of the elements of the basis of the
linear space; by requirement that multiplication is distributive the
definition extends uniquely to general vectors.

For example, for
\begin{multline}
\label{eq:przykladmnozenie} \rho_1=\{1,2,7,8\},\
\rho_2=\{3,4,5,6\},\\  \pi^1=\big\{ \{1,7\}, \{2\},\{8\} \big\},\
\pi^2=\big\{ \{3,5\},\{4,6\} \big\}
\end{multline}
we have
\begin{multline*}\pi^1 \cdot \pi^2 =
\big\{ \{1,7\},\{2,4,6\},\{3,5\},\{8\} \big\}+\\
\big\{ \{1,7\},\{2,4,6\},\{3,5,8\} \big\}+ \big\{
\{1,3,5,7\},\{2,4,6\},\{8\} \big\}.
\end{multline*}

\subsubsection{Geometric interpretation of partitions multiplication}
\label{subsec:geometricmultiplication} 
\begin{figure}[bt]
\includegraphics{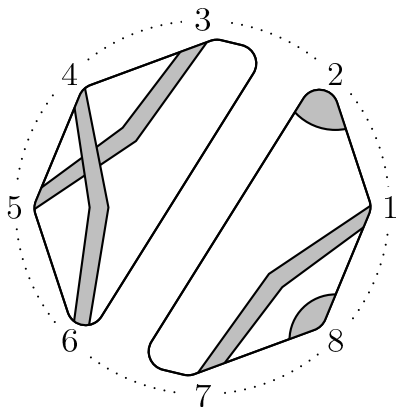}
\caption{Graphical representation of example
\eqref{eq:przykladmnozenie}.} \label{fig:mnozeniepartycji}
\end{figure}

Suppose for simplicity that
$X=\{1,2,\dots,n\}$. On the surface of a large sphere we draw a
small circle on which we mark counterclockwise
 points $1,2,\dots,n$.
Inside the circle we cut $r$ holes; for any $1\leq s\leq r$ the
corresponding hole has a shape of a disc, the boundary of which
passes through the points from the block $\rho_s$ (it is possible to
perform this operation since $\rho$ is non--crossing). For every
$1\leq s\leq r$ the partition $\pi^s$ connects some points on the
boundary of the hole $\rho_s$ and this situation corresponds exactly
to the case we considered in Section \ref{subsec:zaklejanie}. We
shall glue to the hole $\rho_s$ only the first collection of discs
that we considered in Section \ref{subsec:zaklejanie}, i.e.~the
discs which correspond to the blocks of the partition $\pi^s$. Thus
we obtained a number of holes with a collection of glued discs (cf
Figure \ref{fig:mnozeniepartycji}).

\begin{figure}[bt]
\includegraphics{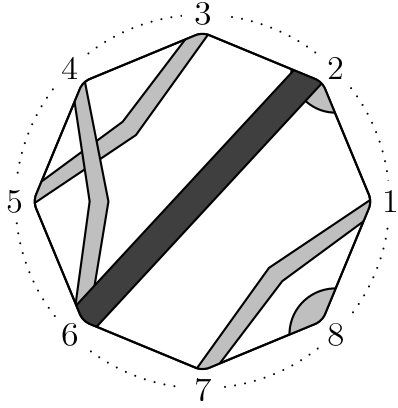}
\caption{Figure \ref{fig:mnozeniepartycji} after inflating small
holes.} \label{fig:mnozeniepartycji2}
\end{figure}

When we inflate the original small holes inside the circle we may
think about this picture alternatively: instead of $r$ small holes
we have a big one (in the shape of the circle) but some arcs on its
boundary are glued by extra discs (on Figure
\ref{fig:mnozeniepartycji2} drawn in black) given by the partition
$\rho_{\compinv}\in\NC(1,2,\dots,n)$. Furthermore we still have all
discs (on Figure \ref{fig:mnozeniepartycji2} drawn in gray)
corresponding to partitions $\pi^s$. We merge discs from these two
collections if they touch the same vertex. After this merging the
collection of discs corresponds to the partition $\rho_{\compinv}
\vee (\pi^1 \cup \cdots \cup\pi^t)$, i.e.~the smallest partition
which is bigger than both $\rho_{\compinv}$ and $\pi^1\cup\cdots\cup
\pi^t$.

The last step is to consider all ways of merging of the discs (or
equivalently: all partitions $\sigma\geq \big( \rho_{\compinv} \vee
(\pi^1 \cup \cdots \cup\pi^t)\big) $) with the property that any two
vertices that were lying on the boundary of the same small hole
$\rho_s$ if were not connected by a disc from the collection $\pi^s$
then they also cannot be connected after all mergings.

\subsection{Map $\Sigma$ is a homomorphism}
Since the map $\Sigma=\fppsinf$ was defined on the basis of the
linear space of partitions, we can extend it linearly in a unique
way. It turns out that $\fppsinf$ preserves also the multiplication
hence it is a homomorphism. However, since we consider only
`$\rho$--ordered' products of partitions, we have to state it in a
bit unusual way.
\begin{claim}
\label{claim:mnozenie} Let $\rho=\{\rho_1,\dots,\rho_r\}$ be a
non--crossing partition of the set $\{1,2,\dots,n\}$ and for every
$1\leq s\leq r$ let $\pi^s$ be a partition of the set $\rho_s$. Then
$$\Sigma\Big( \prod_{s} \pi^s \Big)=\prod_{1\leq s\leq r} \Sigma_{\pi^s},$$
where the multiplication on the left hand side should be understood
as the $\rho$--ordered product of partitions and on the right hand
side it should be understood as the usual product of commuting
elements in $\C(\PartPerm{\infty})$.
\end{claim}
This Claim will follow from Theorem
\ref{theo:sigmajesthomomorfizmem}. Due to distributivity of
multiplication the above Claim remains true if we allow $\pi^s$ to
be a linear combination of partitions. It would be very interesting
to prove Claim \ref{claim:mnozenie} directly if one treats Claim
\ref{claim:definicjasigma} as a definition of $\Sigma_{\pi}$.

The following result was proved in
\cite{KerovOlshanski1994,IvanovKerov1999,LascouxThibon2001vertex,BianeCharacters}.
\begin{corollary}
\label{cor:uniwersalnestale} For any choice of tuples of positive
integers $(k_1,\dots,k_m)$ and $(k'_1,\dots,k'_n)$ there exists a
function $f$ with a finite support defined on the set of the tuples
of integers such that \begin{equation} \label{eq:glod}
\Sigma_{k_1,\dots,k_m} \cdot \Sigma_{k'_1,\dots,k'_n} =
\sum_{p_1,\dots,p_r} f_{p_1,\dots,p_r}
\Sigma_{p_1,\dots,p_r}.\end{equation}
%
\end{corollary}
\begin{proof}
It is enough to find partitions $\pi, \pi'$ such that
$\Sigma_{\pi}=\Sigma_{k_1,\dots,k_m}$,
$\Sigma_{\pi'}=\Sigma_{k'_1,\dots,k'_n}$ and use the fact that
$\Sigma_{\pi} \Sigma_{\pi'}=\Sigma_{\pi \cdot \pi'}$.
\end{proof}

More detailed information about the product $\Sigma_{k_1,\dots,k_m}
\cdot \Sigma_{k'_1,\dots,k'_n}$ is provided by the following result
which was proved by Ivanov and Olshanski \cite{IvanovOlshanski2002}.
\begin{corollary}
\label{cor:filtration} For any choice of integers
$k_1,\dots,k_m,k'_1,\dots,k'_n\geq 1$ let us express the product
$\Sigma_{k_1,\dots,k_m} \cdot
\Sigma_{k'_1,\dots,k'_n}\in\C(\PartPerm{\infty})$ as a linear
combination of normalized conjugacy class indicators $\Sigma$. Then
\begin{multline*}
\Sigma_{k_1,\dots,k_m} \cdot \Sigma_{k'_1,\dots,k'_n} =
\Sigma_{k_1,\dots,k_m,k'_1,\dots,k'_n} + \\
\text{(terms of degree at most
$k_1+\cdots+k_m+k'_1+\cdots+k'_n+m+n-2$)}
\end{multline*}
and therefore \eqref{eq:filtration} indeed defines a filtration.
\end{corollary}


\begin{proof}
Define $l_i=k_i+1$ and $l_i'=k_i'+1$. We set $\pi^1$ to be a
partition of the set of consecutive integers
$\rho_1=\{1,\dots,l_1+\cdots+l_m\}$ which has only one non--trivial
block $\pi^1_1=\{l_1,l_1+l_2,l_1+l_2+l_3,\dots,l_1+\cdots+l_m\}$ and
let $\pi^2$ be a partition of the set of consecutive integers
$\rho_2=\{l_1+\cdots+l_m+1,l_1+\cdots+l_m+2,\dots,l_1+\cdots+l_m+l'_1+\cdots+l'_n\}$
which has only one non--trivial block
$\pi^2_1=\{l_1+\cdots+l_m+l'_1,l_1+\cdots+l_m+l'_1+l'_2,\dots,
l_1+\cdots+l_m+l'_1+\cdots+l'_n\}$; all the other blocks of these
partitions consist of single elements.
One can easily check that $\Sigma_{\pi^1}=\Sigma_{k_1,\dots,k_m}$
and $\Sigma_{\pi^2}= \Sigma_{k'_1,\dots,k'_n}$. We set
$\rho=\{\rho_1,\rho_2\}$; let us compute the $\rho$--ordered product
$\pi^1 \cdot \pi^2$.

Equation \eqref{eq:genusowatosc} implies that the terms of the
maximal degree $l_1+\cdots+l_m+l'_1+\cdots+l'_n$ will correspond to
non--crossing partitions (genus equal to zero) and the degree
corresponding to any crossing partition cannot exceed
$l_1+\cdots+l_m+l'_1+\cdots+l'_n-2$. It is enough to show that there
is only one non--crossing partition $\sigma$ which contributes to
$\pi^1 \cdot \pi^2$, namely  the partition $\tau=(\pi^1 \cup \pi^2)
\vee \rho_{\comp^{-1}}$ which has only one non--trivial block
$\pi^1_1\cup\pi^2_1$. This statement is pretty obvious from the
geometric interpretation; we provide a more detailed proof below.

Let us consider a partition $\sigma\neq \tau$ which contributes to
$\pi^1 \cdot \pi^2$. Of course we have $\sigma\geq\tau$. Let $a,c$
be a pair of elements which are connected by $\sigma$ and are not
connected by $\tau$; by the definition of the product of partitions
$a$ and $c$ cannot be both elements of $\rho_1$ or both elements of
$\rho_2$ so let us assume that $a\in\rho_1$ and $c\in\rho_2$.

Suppose that $a\in\pi^1_1$; then $c\notin\pi^2_1$ (otherwise $a$ and
$c$ would be connected by $\tau$). It follows that there is an
element of $\rho_2\setminus\pi^2_1$, namely $c$, which is connected
by $\sigma$ with $a\in\pi_1^1$, and the latter element must be
connected by $\sigma$ with the elements of $\pi^2_1$. This
contradicts the definition of $\pi^1 \cdot \pi^2$. Therefore
$a\notin\pi^1_1$ and similarly we show that $c\notin\pi^2_1$.

It follows that a tuple $a,b,c,d$ with $b=l_1+\cdots+l_m$,
$d=l_1+\cdots+l_m+l'_1+\cdots+l'_n$ is the one required for $\sigma$
to be crossing.

\end{proof}

\section{Free cumulants of partitions}
\label{sec:mainresult}

\subsection{Free cumulants}
\label{subsec:freeprobability} \index{free cumulants}


\index{cumulants, free} Let us fix some finite ordered set $X$, a
linear space $V$ and a map $M:\NC(X)\rightarrow V$ called the moment
map. Usually, the space $V$ carries some kind of multiplicative
structure and there exists a sequence $M_1,M_2,\dots$, called moment
sequence, such that for every non--crossing partition
$\rho=\{\rho_1,\dots,\rho_r\}$ the moment map is given by a
multiplicative extension
\begin{equation}
\label{eq:momentmappa} M_{\rho}=\prod_{1\leq s\leq r} M_{|\rho_s|}.
\end{equation}

For every integer $n\geq 1$ the Speicher's free cumulant $R_n$
\cite{Speicher1994,Speicher1997,Speicher1998} is defined by
\index{$Moeb$@$\Moeb_{\pi}$} \index{M\"obius function}
\begin{equation}
\label{eq:definicjawolnejkumulanty}
R_{n}=\sum_{\rho\in\NC(1,2,\dots,n)} \Moeb_{\rho_{\comp}} M_{\rho},
\end{equation}
where the sum runs over all non--crossing partitions of the set
$\{1,2,\dots,n\}$ and $\Moeb$ is the M\"obius function on the
lattice of non--crossing partitions given explicitly by
$$\Moeb_{\sigma}=\prod_{1\leq s\leq r} (-1)^{|\sigma_s|-1}
 c_{|\sigma_s|-1}$$
for $\sigma=\{\sigma_1,\dots,\sigma_r\}$ and where
$c_k=\frac{(2k)!}{(k+1)! k!}$ denotes the Catalan number.

We are going to use in this article the following defining property
of the M\"obius function \cite{Speicher1994}.
\begin{lemma}
\label{lem:pocomoebius} For every $\rho\in\NC(1,2,\dots,n)$
$$\sum_{\substack{\pi\in\NC(1,\dots,n)\\ \pi\leq \rho}} \Moeb_{\pi}=
 \begin{cases} 1 & \text{if } \rho \text{ is trivial,} \\
0 & \text{otherwise.} \end{cases}
$$
\end{lemma}

One can show \cite{Speicher1994} that if \eqref{eq:momentmappa}
holds then the moments can be computed from the corresponding
cumulants by the following simple formula
\begin{equation}
\label{eq:momentkumulant} M_n =\sum_{\rho\in\NC(1,\dots,n)}
R_{\rho},
\end{equation}
where similarly as in \eqref{eq:momentmappa} for every non--crossing
partition $\rho=\{\rho_1,\dots,\rho_r\}$ we set
$$
R_{\rho}=\prod_{1\leq s\leq r} R_{|\rho_s|}.
$$


In a typical application the moments $M_n$ are real numbers given by
$$ M_n= \int_{-\infty}^{\infty} x^n \ d\mu(x),$$
where $\mu$ is a probability measure on $\R$ with all moments finite
or given by $M_n=\E(X^n)$, where $X$ is some random variable. It is
easy to check that in this case the free cumulant $R_n$ behaves like
a homogeneous polynomial of degree $n$ in a sense that the free
cumulant $R_n(D^p \mu)$ of the dilated measure is related to the
free cumulant $R_n(\mu)$ of the original measure by
$$ R_n(D^p \mu)= p^n R_n(\mu).$$
This very simple scaling property is very useful in the study of
asymptotic properties. Another advantage of free cumulants is that
the generating function (called Voiculescu's $R$--transform
\cite{VoiculescuDykemaNica}) of the sequence of free cumulants of a
given measure $\mu$ can be computed directly from the Cauchy
transform of $\mu$ which is very useful in practical applications.

The definition of the free cumulants itself does not involve the
multiplicative structure of $V$; this structure is accessed only
implicitly through the moment map $M$. Therefore a question arises
if there are some reasonable generalizations of
\eqref{eq:momentmappa}. The key point is that we need to define the
moments $M_\rho$ only for non--crossing partitions $\rho$ and it is
well--known that non--crossing partitions of $n$--element set are in
correspondence with certain ways of writing brackets in a product of
$n$ factors and hence recursive definitions of $M_\rho$ can be
successfully applied. Canonical examples are provided by
operator--valued free probability \cite{Speicher1998}. In the
following we will use this weakness of requirements for the moment
map $M$ in order to define and study the partition--valued free
cumulants.

\subsection[Free cumulants of partitions]%
{Free cumulants of partitions and of the Jucys--Murphy element}

\subsubsection{Moments and cumulants of partitions} \label{sec:momentypartycji} Let
$X$ be an ordered set. We define the moment $\Mpp_X$ (where the
letter $\pp$ stands for partition) by \index{$MPP$@$\Mpp$}
$$\Mpp_X=\sum_{\sigma\in\PP(X)} \sigma,$$
where the sum denotes a formal linear combination and it runs over
all partitions of $X$.

Let $\rho=\{\rho_1,\dots,\rho_r\}$ be a non--crossing partition.
Motivated by \eqref{eq:momentmap} we define
$$\Mpp_\rho=\prod_{s}  \Mpp_{\rho_s}=\sum_{\sigma^1,\dots,\sigma^r} \prod_s \sigma^s,$$
where the sum runs over all tuples $(\sigma^1,\dots,\sigma^r)$ such
that $\sigma^s$ is a partition of the set $\rho_s$ and products
denote $\rho$--ordered products of partitions.

\begin{proposition}
\label{prop:formulanamomentypartycji} For every non--crossing
partition $\rho$ of the set $X$ we have
\begin{equation}
\label{eq:momentpartition} \Mpp_\rho=\sum_{\sigma\geq
\rho_{\compinv}} \sigma,
\end{equation} where the sum runs over all partitions $\sigma$ of
the set $X$ such that $\sigma\geq \rho_{\compinv}$.
\end{proposition}
\begin{proof}
Probably the best way to prove the proposition is to consider the
graphical description of the partition multiplication from Section
\ref{subsec:geometricmultiplication}. We present an alternative
proof below.

 From the very definition of multiplication one can see that
every partition $\sigma$ which appears in $\Mpp_\rho$ must fulfill
$\sigma\geq \rho_{\compinv}$.

To any partition $\sigma$ of the set $X$ we can assign a tuple
$(\sigma^1,\dots,\sigma^r)$, where $\sigma^s$ (a partition of the
set $\rho_s$) is obtained from $\sigma$ by restriction to the set
$\rho_s$. Now it is easy to observe that every $\sigma\geq
\rho_{\compinv}$ appears exactly once in $\prod_s \Mpp_{\rho_s}$,
namely in the factor $\prod_s \sigma^s$.
\end{proof}

The moment map $\Mpp$ gives rise to a sequence of free cumulants
which will be denoted by $\Rpp_n$.

\subsubsection{Moments and cumulants of Jucys--Murphy
element}

\label{subsec:momentmap}
%

Inspired by \eqref{eq:momentmappa} we shall denote by $\MJM_\rho$
the partition--indexed moments of the Jucys--Murphy element, defined
by
\begin{equation}
\label{eq:momentmap} \MJM_{\rho}=\prod_{1\leq s\leq r}
\MJM_{|\rho_s|}=\prod_{1\leq s\leq r}
\E(J^{|\rho_s|})\in\C(\PartPerm{\infty}),
\end{equation}
where $\rho=\{\rho_1,\dots,\rho_r\}$ is a non--crossing partition.
The moment map $\MJM$ gives rise to a sequence of free cumulants
which will be denoted by $\RJM_n$. Please note that $\MJM_1=\E(J)=0$
and therefore also $\RJM_1=0$.


\subsubsection{The relation between the partitions and the Jucys--Murphy
element} The following Claim provides a link between the moments of
the Jucys--Murphy element and the moments of partitions.
\begin{claim}
\label{claim:dobryjurcysmurphy} Let $X$ be a finite ordered set, let
$\rho$ be a non--crossing partition of a finite ordered set and
$n\geq 1$ be an integer.
Then \begin{align} \label{eq:claim1} \Sigma(\Mpp_{X}) &=\MJM_{|X|}, \\
\label{eq:claim2} \Sigma(\Mpp_\rho)&=\MJM_\rho, \\
\label{eq:claim3} \Sigma(\Rpp_n)&=\RJM_n, \end{align} where $\Sigma$
is the map considered in Section
\ref{subsec:pierwszawzmiankaofppsq}.
\end{claim}

We postpone the proof of \eqref{eq:claim1} to Theorem
\ref{theo:zgodnoscpartycjizjurcysiem}. Equation \eqref{eq:claim1}
implies \eqref{eq:claim2} immediately because  $\Sigma$ is a
homomorphism. Equation \eqref{eq:claim3} is an immediate consequence
of \eqref{eq:claim2}.

This result implies that all questions about the moments of the
Jucys--Murphy element can be easily translated into questions about
moments of partitions and the machinery of Section
\ref{sec:calculus} can be applied.

\subsubsection{Free cumulants of partitions}


\begin{proposition}  For every $n\geq 1$
\begin{equation}
\label{eq:kumulantypartycji} \Rpp_n=\sum_{\pi\in\PP(1,2,\dots,n)}
I_{\pi} \  \pi,
\end{equation}
\begin{equation}
\label{eq:kumulantypartycji2} \RJM_n=\sum_{\pi\in\PP(1,2,\dots,n)}
I_{\pi} \  \Sigma_\pi,
\end{equation}
where $I_{\pi}\in\Z$, called\index{free index of a
partition}\index{partition!free index of} \index{$I_\pi$} free index
of $\pi$, is defined by
\begin{equation}
\label{eq:formulanakrotnosc}
I_{\pi}=\sum_{\substack{\rho\in\NC(1,2,\dots,n) \\ \rho\leq \pi}}
\Moeb_{\rho}.
\end{equation}
\end{proposition}
\begin{proof} From \eqref{eq:definicjawolnejkumulanty} and
\eqref{eq:momentpartition} it follows that
\begin{multline*}
\Rpp_n=\sum_{\rho\in\NC(1,2,\dots,n)} \Moeb_{\rho_{\comp}}
\sum_{\substack{\pi\in\PP(1,2,\dots,n)\\ \pi\geq \rho_{\comp}} }
\pi=\\
\sum_{\pi\in\PP(1,2,\dots,n)} \pi
\sum_{\substack{\rho\in\NC(1,2,\dots,n) \\ \rho_{\comp} \leq \pi}}
\Moeb_{\rho_{\comp}}= \sum_{\pi\in\PP(1,2,\dots,n)} \pi
\sum_{\substack{\rho\in\NC(1,2,\dots,n) \\ \rho \leq \pi}}
\Moeb_{\rho},
\end{multline*}
where in the last equality we used the fact that
$\rho\mapsto\rho_{\comp}$ is a permutation of the set of
non--crossing partitions.

The second identity \eqref{eq:kumulantypartycji2} follows now from
Claim \ref{claim:dobryjurcysmurphy}. \end{proof}


\subsubsection{Evercrossing partitions} \index{evercrossing partition}
\index{partition!evercrossing}

\begin{figure}[bt]
\includegraphics{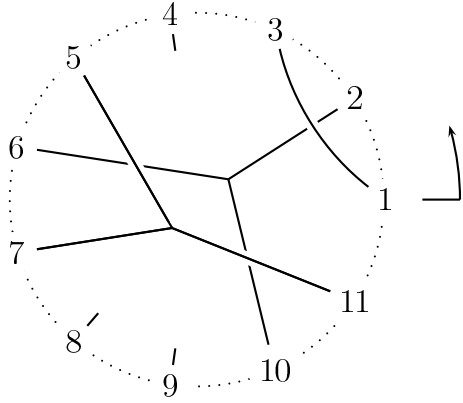}
\caption{Example of an evercrossing partition}
\label{fig:evercrossing}
\end{figure}

We say that a partition $\pi$ of a finite ordered set $X$ is
evercrossing if for any $a<c$ ($a,c\in X$) such that $a$ is
connected with $c$ by $\pi$ there exist $b,d\in X$ such that
\begin{enumerate}
\item elements $a,b,c,d$ are ordered up to a cyclic rotation,
i.e.~either $a<b<c<d$ or $d<a<b<c$, \item elements $b$ and $d$ are
connected by $\pi$, \item elements $a,b,c,d$ are not elements of the
same block of $\pi$.
\end{enumerate}
Up to some technical details, this means that on the graphical
representation of $\pi$ every line connecting a pair of elements
must be crossed by some other line. For example, partition from
Figure \ref{fig:evercrossing} is evercrossing, but partition from
Figure \ref{fig:przecinajacapartycja} is not.

The following result shows that only evercrossing partitions
contribute to sums \eqref{eq:kumulantypartycji} and
\eqref{eq:kumulantypartycji2}.
\begin{theorem}
\label{theo:wiwatevercrossing} If a partition $\pi$ is not
evercrossing then $I_{\pi}=0$.
\end{theorem}

In order to prove this result we will need the following lemma.

\begin{lemma}
\label{lem:fundamentalny} Let $\pi=\{\pi_1,\dots,\pi_r\}$ be a
partition of $X$. Let $\sigma=\{\pi_1,\dots,\pi_s\}$,
$\tau=\{\pi_{s+1},\dots,\pi_r\}$ be a decomposition of $\pi$ into
two partitions of disjoint sets $S=\pi_1\cup\cdots\cup\pi_s$ and
$T=\pi_{s+1}\cup\cdots\cup\pi_{r}$ such that $\sigma$ is
non--crossing.

For $\rho\in\NC(T)$ we shall denote (with a small abuse of notation)
by $\rho_{\comp}\in\NC(S)$ the biggest non--crossing partition with
a property that $\rho\cup\rho_{\comp}$ is non--crossing. We denote
by $\sigma\wedge\rho_{\comp}\in\NC(S)$ the biggest non--crossing
partition which is smaller both than $\sigma$ and $\rho_{\comp}$.

Then the free index of $\pi$ is given by
\begin{equation}
\label{eq:uproszczeniekumulanty}
I_{\pi}=\sum_{\substack{\rho\in\NC(T)\\ \rho\leq \tau
\\ (\sigma \wedge \rho_{\comp}) \text{ is trivial}  }}
  \Moeb_{\rho}.
\end{equation}
\end{lemma}
\begin{proof}
Since every $\rho\in\NC(X)$, $\rho\leq \pi$  can be decomposed as
$\rho=\tilde{\rho}\cup\hat{\rho}$, where $\tilde{\rho}\in\NC(S)$,
$\tilde{\rho}\leq \sigma$ and $\hat{\rho}\in\NC(T)$, $\hat{\rho}\leq
\tau$ it follows from \eqref{eq:formulanakrotnosc} that
$$ I_{\pi}=\sum_{\substack{\hat{\rho}\in \NC(T)\\ \hat{\rho}\leq \tau}}  \Moeb_{\hat{\rho}}
\Big( \sum_{\substack{  \tilde{\rho}\in \NC(S) \\ \tilde{\rho} \leq
(\sigma\wedge \hat{\rho}_{\comp}) }} \Moeb_{\tilde{\rho}} \Big).
$$
We apply Lemma \ref{lem:pocomoebius} and observe that the second sum
is equal to zero unless $\sigma\wedge\hat{\rho}_{\comp}$ is trivial;
otherwise it is equal to $1$.
\end{proof}

\begin{proof}[Proof of Theorem \ref{theo:wiwatevercrossing}]
Let  $a,c$ be the vertices which have the property required for
$\pi=\{\pi_1,\dots,\pi_r\}$ not to be evercrossing. By a change of
numbering of blocks we can always assume that $a,c\in\pi_1$.
Following the notation of Lemma \ref{lem:fundamentalny} we set
$S=\pi_1$, $T=\pi_2\cup\cdots\cup\pi_r$, $\sigma=\{\pi_1\}$,
$\tau=\{\pi_2,\dots,\pi_r\}$.

Observe that for every ${\rho}\in\NC(T)$ such that ${\rho}\leq
\tau$, we have that the partition ${\rho} \cup \big\{ \{a,c\}
\big\}$ is non--crossing, therefore ${\rho}_{\comp}\geq
\big\{\{a,c\} \big\}$ and  ${\rho}_{\comp} \wedge \sigma \geq
\big\{\{a,c\} \big\}$ must contain a non--trivial block and the sum
in \eqref{eq:uproszczeniekumulanty} contains no summands.
\end{proof}
\begin{remark}
The converse implication is not true, as one can see on an example
from Figure \ref{fig:kontrprzyklad}.

\begin{figure}[bt]
\includegraphics{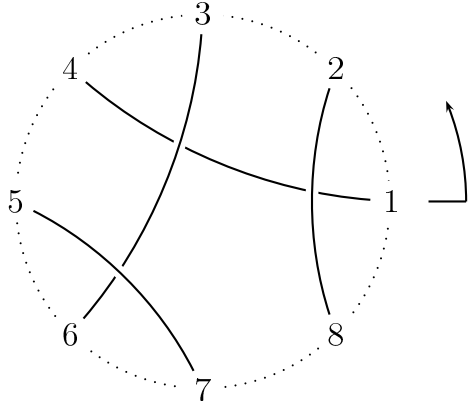}
\caption{Example of an evercrossing partition with the free index
equal to $0$.} \label{fig:kontrprzyklad}
\end{figure}

\end{remark}

\subsection{First--order asymptotics of free cumulants}

\begin{theorem}
\label{theo:pierwszystopien}
Let integers $k_1,\dots,k_t\geq 2$ be given and let us express
$\RJM_{k_1} \cdots \RJM_{k_t}\in\myalgebra$ as a linear combination
of normalized conjugacy class indicators $\Sigma$. Then
\begin{multline}
\label{eq:highestorde} \RJM_{k_1} \cdots
\RJM_{k_t}=\Sigma_{k_1-1,\dots,k_t-1}+\\ \text{(terms of degree at
most } k_1+\cdots+k_t-2\text{)},
\end{multline}
where the degree is taken with respect to the filtration
\eqref{eq:filtration}.
\end{theorem}
\begin{proof}
Let us consider the case $t=1$; we set $n=k_1$. Propositions
\ref{prop:genusowatosc} and \ref{prop:genuszero} imply that the
highest--order terms in the expansion \eqref{eq:kumulantypartycji2}
correspond to non--crossing partitions and that the degree of all
remaining terms will be at most $n-2$; by Theorem
\ref{theo:wiwatevercrossing} a partition contributes in the sum
\eqref{eq:kumulantypartycji2} only if it is evercrossing. In order
to find the leading term we need to find all non--crossing
partitions of the set $\{1,2,\dots,n\}$ which are at the same time
evercrossing. This combination of adjectives sounds oxymoronic which
suggests that there should not be too many of such partitions and
indeed it is easy to check that the only such partition is the
trivial one. Hence
\begin{multline*}
\RJM_n=\Sigma_{ \{ \{1\},\{2\},\dots,\{n\} \} }+\text{(terms of
degree at most } n-2\text{)} =
\\ \Sigma_{n-1} + \text{(terms of degree at most } n-2\text{)}
\end{multline*}
which finishes the proof of the case $t=1$.

The general case follows easily from Corollary \ref{cor:filtration}.
\end{proof}

\begin{corollary}
\label{cor:algebraicallyfree} Elements $(\RJM_k)_{k\geq 2}$ are
algebraically free; also elements $(\MJM_k)_{k\geq 2}$ are
algebraically free and therefore the gradation \eqref{eq:gradation}
is well--defined.
\end{corollary}
\begin{proof}
Let $P(x_2,x_3,\dots)\neq 0$ be a polynomial; our goal is to show
that $P(\RJM_2,\RJM_3,\dots)\in\C(\PartPerm{\infty})$ is non--zero.
In order to show this we express $P(\RJM_2,\RJM_3,\dots)$ as a
linear combination of conjugacy class indicators $\Sigma$ and
Theorem \ref{theo:pierwszystopien} shows that the highest--order
terms do not cancel.

In order to show that $(\MJM_k)_{k\geq 2}$ are algebraically free we
consider a polynomial $Q(x_2,x_3,\dots)\neq 0$. We apply
\eqref{eq:momentkumulant} and obtain a polynomial $P(x_2,x_3,\dots)$
such that $Q(\MJM_2,\MJM_3,\dots)=P(\RJM_2,\RJM_3,\dots)$. Equation
\eqref{eq:definicjawolnejkumulanty} implies that
$P(x_2,x_3,\dots)\neq 0$. It follows that
$Q(\MJM_2,\MJM_3,\dots)=P(\RJM_2,\RJM_3,\dots)\neq 0$.
\end{proof}

The following result was proved in
\cite{Biane1998,IvanovOlshanski2002}.
\begin{theorem}
\label{theo:kontrolastopnia} For every tuple $k_1,\dots,k_t$ of
positive integers there exists a polynomial
$K_{k_1,\dots,k_t}(\RJM_2,\RJM_3,\dots)$, called \index{Kerov
polynomial} Kerov polynomial, such that
$$ \Sigma_{k_1,\dots,k_t}=K_{k_1,\dots,k_t}(\RJM_2,\RJM_3,\dots). $$
In other words, each of the families $(\RJM_n)_{n\geq 2}$,
$(\MJM_n)_{n\geq 2}$ and
$(\Sigma_{k_1,\dots,k_m})_{k_1,\dots,k_m\geq 1}$ generates the same
commutative algebra $\myalgebra$.

With respect to the gradation \eqref{eq:gradation} the leading term
of the Kerov polynomial is given by
\begin{multline}
\label{eq:eksplozjamozgu}
\Sigma_{k_1,\dots,k_t}=\RJM_{k_1+1} \cdots \RJM_{k_t+1} +\\
\text{(terms of degree at most $(k_1+1)+\cdots+(k_t+1)-2$)}
\end{multline}
and
$$\deg \Sigma_{k_1,\dots,k_t}=(k_1+1)+\cdots+(k_t+1).$$
In other words, the filtration \eqref{eq:filtration} is induced by
gradation \eqref{eq:gradation}.
\end{theorem}
\begin{proof}
We shall prove this theorem by induction over
$(k_1+1)+\cdots+(k_t+1)$.


Let us express $\Sigma_{k_1,\dots,k_t}-\RJM_{k_1+1} \cdots
\RJM_{k_t+1}$ as a linear combination of normalized conjugacy class
indicators $\Sigma$. Theorem \ref{theo:pierwszystopien} implies that
with respect to the filtration \eqref{eq:filtration} all these terms
have degree at most $(k_1+1)+\cdots+(k_t+1)-2$ and hence the
inductive hypothesis can be applied: every of these terms can be
expressed by the corresponding Kerov polynomial. Observe that from
the very definition of free cumulants
\eqref{eq:definicjawolnejkumulanty} it follows that $\deg
\RJM_{k_1+1} \cdots \RJM_{k_t+1}=(k_1+1)+\cdots+(k_t+1)$ with
respect to the gradation \eqref{eq:gradation} and hence the theorem
follows.
\end{proof}

\subsection{Higher order expansion for free cumulants}
\label{sec:evercrossing}

In this section we present a method of enumerating all evercrossing
partitions with a given genus. The main concept is to simplify a
given evercrossing partition by a number of steps. After these
simplifications we obtain an evercrossing pair partition (with some
extra properties) which can be enumerated easily. By reversing the
process of simplifications we shall therefore obtain all
evercrossing partitions with a given genus.

\subsubsection{Simplification, step 1} \label{subsec:simple1} 
\begin{figure}[bt]
\includegraphics{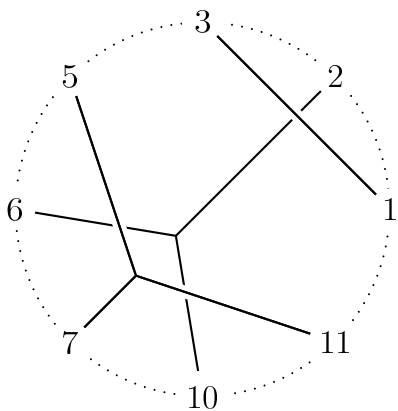}
\caption[Partition from Figure \ref{fig:evercrossing} after the
first step of simplification]{Partition from Figure
\ref{fig:evercrossing} after the first step of the simplification
algorithm.} \label{fig:wynikstep1}
\end{figure}

Let an
evercrossing partition be given. In the first step we remove all its
trivial blocks. For example, Figure \ref{fig:wynikstep1} depicts the
outcome of this step performed on an evercrossing partition from
Figure \ref{fig:evercrossing}. In fact we do not really care what
are the labels of vertices and it would be a good idea not to write
them down at all; in this example we do write them down for the
reason of clarity.

\subsubsection{Simplification, step 2} 
\begin{figure}[bt]
\includegraphics{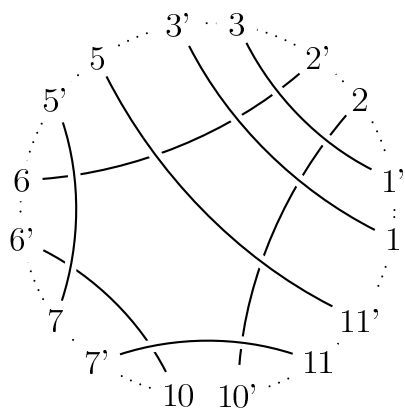}
\caption[Partition from Figure \ref{fig:evercrossing} after the
second step of the simplification]{Partition from Figure
\ref{fig:evercrossing} after the second step of the simplification
algorithm.} \label{fig:wynikstep2}
\end{figure}

\begin{figure}[bt]
\includegraphics{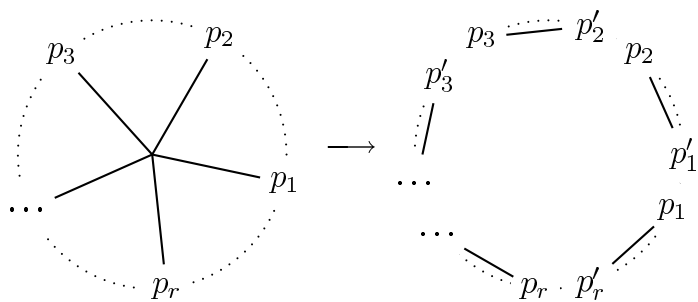}
\caption[Elementary operation of step 2 of the simplification
algorithm]{Elementary operation of step 2 of the simplification
algorithm: one of the blocks is replaced by its fat version.}
\label{fig:step2elementary}
\end{figure}

In the second step of
simplification, we consider the fat partition corresponding to the
outcome of the first step (cf example on Figure
\ref{fig:wynikstep2}). In other words, we perform the operation
depicted on Figure \ref{fig:step2elementary} on each of the blocks.

\subsubsection{Simplification, step 3} \label{subsec:simple3} 
\begin{figure}[bt]
\includegraphics{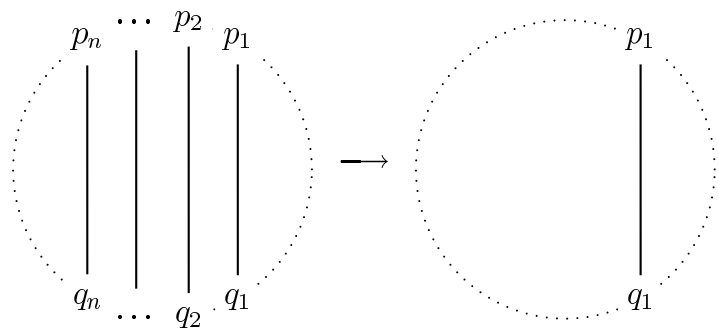}
\caption{The elementary operation of step 3 of the simplification
algorithm: a number of parallel lines is replaced by a single one.}
\label{fig:step3}
\end{figure}

\begin{figure}[bt]
\includegraphics{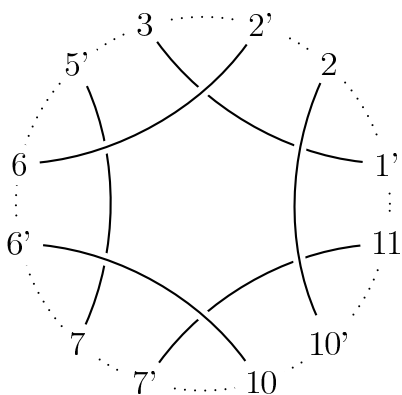}
\caption{The outcome of the step 3 of the simplifying algorithm for
partition from Figure \ref{fig:evercrossing}.}
\label{fig:wynikstep3}
\end{figure}

The
outcome of the second step is a pair partition. In the third step of
the simplification algorithm we perform the operation depicted on
Figure \ref{fig:step3}. To be precise: for $n\geq 2$ let
$(p_1,\dots,p_n)$ and $(q_n,q_{n-1},\dots,q_1)$ be two sequences of
consecutive vertices of this pair partition such that $p_k$ is
connected with $q_k$.  We remove all vertices
$p_2,p_3,\dots,p_n,q_2,q_3,\dots,q_n$ and leave vertices $p_1$ and
$q_1$ unchanged. We iterate this removal of vertices as long as it
is possible (cf example on Figure \ref{fig:wynikstep3}).

A careful reader might observe that after rotating the Figure
\ref{fig:step3} upside--down the roles played by $p_1$, $q_1$ are
interchanged with $q_n$ and $p_n$ therefore the above rule is not
very precise which labels should be carried by the vertices of the
surviving block; but in fact we do not really care about these
labels. In particular, we do not care which of the labels is the
smallest, therefore we identify pair partitions which are cyclic
rotations of each other.

\subsubsection{Outcome of the simplification algorithm}
\begin{theorem}
\label{theo:gorneszacowanie} The outcome of the simplification
algorithm from Sections \ref{subsec:simple1}--\ref{subsec:simple3}
is an evercrossing pair--partition which has the same genus and free
index as the original partition.

Let $n$ denote the number of vertices of the outcome pair partition
$\pi$. If $\genus_\pi\neq 0$ then
\begin{equation}
n\leq 12 \genus_{\pi}-6. \label{eq:szacowanie}
\end{equation}
\end{theorem}
\begin{proof}
In order to show that all three steps of the simplification
algorithm preserve the genus of a partition, it is enough to notice
that the surface without the boundary we constructed in Section
\ref{subsec:zaklejanie} changes in each of the steps into a
homeomorphic one.

We shall prove now that all three steps preserve the free index of a
partition. Let $\pi$ be a partition of some set $X$, let $x\notin
X$, then $\pi'=\pi\cup \{ x \}$ is a partition of the set $X \cup
\{x\}$ which has an additional trivial block. Observe that every
partition $\rho'\leq \pi'$ must have a form $\rho'=\rho \cup \big\{
\{x\} \big\}$ where $\rho\leq\pi$. Therefore
$$I_{\pi'}=\sum_{\substack{ \rho'\in\NC(X \cup \{x\}) \\ \rho'\leq \pi'}} \Moeb_{\rho'} =
\sum_{\substack{ \rho\in\NC(X) \\ \rho\leq \pi}} \Moeb_{\rho} =
I_{\pi}.$$ It follows that two partitions which differ with one
trivial block have the same free index and, by iterating, the first
step preserves the free index.

The second and the third step can be divided into a number of
elementary operations
(for the second step the elementary operation is the replacement of
only one of the blocks of the original partition by the
corresponding `fat block', cf Figure \ref{fig:step2elementary}, for
the third step the elementary operation is the performance of only
one operation from Figure \ref{fig:step3}). Our goal is to prove
that every elementary operation preserves the free index. Let $\pi$
be a partition and let $\pi'$ be an outcome of one of the elementary
operations performed on $\pi$. We can decompose $\pi=\sigma \cup
\tau$, $\pi'=\sigma'\cup\tau$ (we use the notations of Lemma
\ref{lem:fundamentalny}), where $\tau\in\PP(T)$ is the set of blocks
of $\pi$ not changed by the elementary operation. One can easily
check that for both elementary operations $\sigma\in\NC(S)$ and
$\sigma'\in\NC(S')$ are non--crossing partitions, furthermore for
every $\rho\in\PP(T)$ we have that $(\sigma \wedge
\rho_{\comp})\in\NC(S)$ is trivial if and only if $(\sigma' \wedge
\rho_{\comp})\in\NC(S')$ is trivial. From Lemma
\ref{lem:fundamentalny} applied to $\pi$ and $\pi'$ it follows that
$$ I_{\pi}=\sum_{\substack{\rho\in\NC(T), \\ \rho\leq \tau,
\\ (\sigma \wedge \rho_{\comp}) \text{ is trivial}  }}
  \Moeb_{\rho}=
\sum_{\substack{\rho\in\NC(T), \\ \rho\leq \tau,
\\ (\sigma' \wedge \rho_{\comp}) \text{ is trivial}  }}
  \Moeb_{\rho}= I_{\pi'}$$
and hence both elementary operations preserve the free index.

Observe also, that the first and the third step of the algorithm
preserve clearly the partition property of being evercrossing. It
should be clear from the graphical interpretation that also the
second step preserves this property; nevertheless we provide a more
detailed proof below. Let a partition $\pi$ be an outcome of the
first step of the algorithm. Any pair connected by $\pi_{\fat}$ must
be of the form $(\pi_{s,t}', \pi_{s,t+1})$ (we use the notation from
Section \ref{subsec:fat}). As an outcome of the first step of the
simplification algorithm, $\pi$ contains no trivial blocks hence
$\pi_{s,t}\neq \pi_{s,t+1}$. Since $\pi$ is evercrossing, there
exist $b,d\in\pi_{u}$, $u\neq s$, such that
$\pi_{s,t},b,\pi_{s,t+1},d$ are ordered up to a cyclic rotation. In
other words: on the graphical representation on a circle, elements
of the block $\pi_u$ are not all on the same side of the line which
passes through $\pi_{s,t}$ and $\pi_{s,t+1}$. Therefore there must
be some consequent elements $\pi_{u,v}$, $\pi_{u,v+1}$ of the block
$\pi_u$ such that $\pi_{s,t},\pi_{u,v},\pi_{s,t+1},\pi_{u,v+1}$ are
ordered up to a cyclic rotation. It is easy to see that
$b=\pi_{u,v}'$ and $d=\pi_{s,t+1}$ are the vertices required by the
property of $\pi_{\fat}$ to be evercrossing.

Let $\pi$ be the outcome of the simplification algorithm and such
that $\genus_{\pi}\neq 0$. We keep notation from Section
\ref{subsec:zaklejanie}, i.e.~$r=\frac{n}{2}$ denotes the number of
blocks of $\pi$ and $t$ denotes the number of cycles of
$\pi_{\fat}\circ c$ (or, alternatively, the number of holes after
gluing the first collection of discs). Observe that on the boundary
of every hole there must be at least three intervals which touch
some discs from the first collection (if some hole does not touch
any discs then $\pi=\emptyset$ and $\genus_{\pi}=0$; if some hole
touches only one disc we denote by $a$, $c$ the only two vertices
which are on the boundary of the hole and clearly $(a,c)$ is a pair
of vertices required for $\pi$ to be non--evercrossing; if some hole
touches exactly two discs then it is possible to perform some
simplifications of step 3; a careful reader might observe that it is
possible and perfectly legal for a hole to touch some of the discs
twice, cf the left--hand side partition on Figure
\ref{fig:zwyczajninadzwyczajni}). Since $\pi$ is a pair partition,
every disc of the first collection touches the holes on exactly two
intervals. It follows that $n= 2r \geq 3t $. Equation
\eqref{eq:genus} implies that
$$ n \leq 12 \genus_{\pi}-6$$
which finishes the proof.
\end{proof}

\subsubsection{Interpretation of the free index?} Free index of a
partition is a quite strange combinatorial object. We can think that
it measures how far a given partition is away from non--crossing
partitions or how much it is evercrossing. Its particularly
interesting feature appears in Theorem \ref{theo:gorneszacowanie},
namely that it is being preserved by a number of operations one can
perform on a partition. The list of such operations is far from
being complete and we leave it as an exercise to the reader to find
at least one such operation which does not appear in the
simplification algorithm. Suppose we have completed such a list of
natural operations preserving free index; now we can treat them as
analogues of Reidemeister moves in the knot theory and say that two
partitions are isotopic when one can be obtained from another by a
sequence of such elementary operations. Of course free index will be
an invariant of the isotopy class but we should not expect that it
will always be able to distinguish different isotopy classes. A
question arises: is there some natural (geometric?) interpretation
of such isotopy?

\subsection{The second main result:
the second--order expansion of free cumulants and characters}
\label{sec:secondorderexp}
\subsubsection{Evercrossing partitions with genus $1$}
Theorem \ref{theo:gorneszacowanie} allows us to enumerate (in finite
time) all possible outcomes of the simplifying algorithm for a fixed
genus of the original evercrossing partition. By reversing all three
steps of the simplifying algorithm we are able to write a number of
``templates'' which describe all evercrossing partitions with a
prescribed genus. Let us have a look on an example.

\begin{theorem}
\label{theo:klasyfikacja} Every evercrossing partition with genus
$1$ must be either of the form depicted on Figure \ref{fig:krzyz}
(in this case it has the free index equal to $-1$) or Figure
\ref{fig:gwiazdadawida} (in this case it has the free index equal to
$-2$). On both figures only non--trivial blocks were shown.

\begin{figure}[tb]
\includegraphics{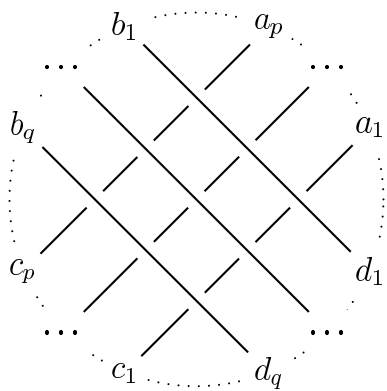}
\caption[Evercrossing partitions with genus $1$ and free index
$-1$]{General pattern of evercrossing partitions with genus $1$ and
free index equal to $-1$ (only non--trivial blocks were shown).}
\label{fig:krzyz}
\end{figure}

\begin{figure}[tb]
\includegraphics{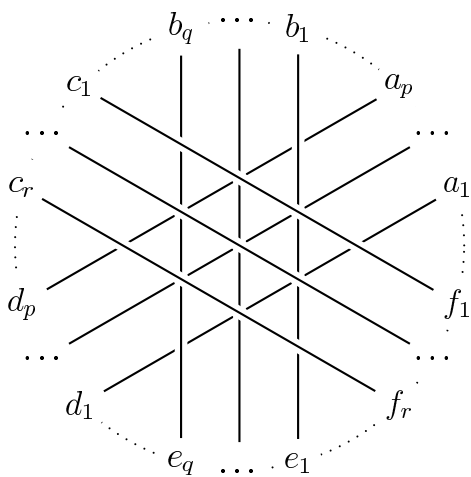}
\caption[Evercrossing partitions with genus $1$ and free index
$-2$]{General pattern of evercrossing partitions with genus $1$ and
free index equal to $-2$ (only non--trivial blocks were shown).}
\label{fig:gwiazdadawida}
\end{figure}

\end{theorem}
\begin{proof}
Theorem \ref{theo:gorneszacowanie} implies that an outcome of the
simplification algorithm obtained for an evercrossing partition with
genus $1$ must have at most $6$ vertices. By direct inspection of
all such pair partitions one can find that there are only two
possible outcomes of the simplifying algorithm and they are depicted
on Figure \ref{fig:zwyczajninadzwyczajni}.

\begin{figure}[bt]
\includegraphics{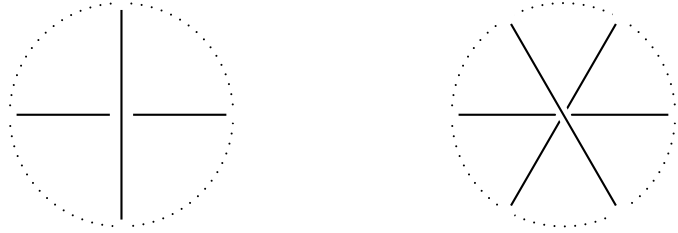}
\caption[Possible outcomes of the simplification algorithm for an
evercrossing partition with genus $1$]{The only two possible
outcomes of the simplification algorithm for an evercrossing
partition with genus $1$. The pair partition on the left--hand side
has free index equal to $-1$ and the pair partition on the
right--hand side has free index equal to $-2$.}
\label{fig:zwyczajninadzwyczajni}
\end{figure}

By reverting the third step of the simplifying algorithm we see that
the outcome of the first two steps of the simplifying algorithm must
be a pair partition either of the form from Figure \ref{fig:krzyz}
or Figure \ref{fig:gwiazdadawida}. For Figure
\ref{fig:gwiazdadawida} we have to study $2^4$ cases now: each of
the numbers $p,q,r$ can be either even or odd, we also have to guess
which of the vertices carry primed and which carry non--primed
labels (we have two possibilities); the case of Figure
\ref{fig:krzyz} is slightly simpler. For each of the $2^4$ cases we
try to reverse the second step of the algorithm (and for most of the
cases it is not possible; a less patient reader will find easily
arguments which would allow to find only the few cases when it is
possible). It turns out that after the first step of the
simplification algorithm must be a pair partition either of the form
depicted on Figure \ref{fig:krzyz} or on Figure
\ref{fig:gwiazdadawida}. Reverting the first step means that we are
allowed to add any number of trivial blocks, which finishes the
proof.
\end{proof}

\subsubsection{Second--order asymptotics of free cumulants}

We leave the proof of the following simple lemma to the reader.
\begin{lemma}
\label{lem:summertime} For any $y,n\in\N$ there are
$\binom{y-1}{n-1}$ solutions of the equation $x_1+\cdots+x_n=y$ if
we require $x_1,\dots,x_n$  to be positive integers.
\end{lemma}

\begin{theorem}
\label{theo:gagggag} For every $n\in\N$ we have
\begin{multline}
\label{eq:dukakao}
 \RJM_{n}=
 \Sigma_{n-1}- \!\!\!\! \sum_{\substack{m_2,m_3,\dots\geq 0 \\
2m_2+3 m_3+4 m_4+\cdots=n-2}} \!\!\!\!
\frac{n (n-1) (n-2)}{24} \binom{m_2+m_3+\cdots}{m_2,m_3,\dots} \times \\
\prod_{s\geq 2} \big( (s-1) \big)^{m_s}
\Sigma_{\underbrace{1,\dots,1}_{m_2 \text{ times}},
\underbrace{2,\dots,2}_{m_3 \text{ times}},
\underbrace{3,\dots,3}_{m_4 \text{ times}},\dots} +\\
\text{(terms of degree at most $n-4$)}
  \end{multline}

\end{theorem}
\begin{proof}
The highest--order term in expansion \eqref{eq:kumulantypartycji2}
was already found in Theorem \ref{theo:pierwszystopien}. The next
highest--order terms correspond to evercrossing partitions with
genus $1$. We shall calculate first the contribution of the
partitions of the form depicted on Figure \ref{fig:krzyz}.

For the partition from Figure \ref{fig:krzyz} the tuple
$(k_1,\dots,k_{p+q-1})$ given by the algorithm from Section
\ref{subsec:pierwszawzmiankaofppsq} is equal to
\begin{multline}
\label{eq:zlepyciagi}
\Big( \big( (a_{s+1}-a_{s})+(c_s-c_{s+1})-1 \big)_{s=1,2,\dots,p-1},\\
\big((b_1-a_p)+(c_p-b_q)+(d_q-c_1)+(a_1-d_1)\big)-3, \\
\big( (b_{s+1}-b_{s})+(d_s-d_{s+1})-1 \big)_{s=1,2,\dots,q-1} \Big).
\end{multline}
Here and in the following, the subtraction is taken modulo $n$, i.e.
$$y-x=\begin{cases} y-x & \text{if } y>x, \\ (y+n)-x & \text{if } y<x. \end{cases}$$
This subtraction has a natural interpretation as the number of
elements between $x$ and $y$ going counterclockwise from $x$ to $y$.

Let $m_2,m_3,\dots$ be a sequence of nonnegative integers such that
$2 m_2+3 m_3+\cdots=n-2$. We shall  count now for how many different
partitions of the form depicted on Figure \ref{fig:krzyz} the tuple
\eqref{eq:zlepyciagi} is equal (up to a permutation) to
\begin{equation}
\label{eq:gromada}
(k_1',k_2',\dots,k_{p+q-1}')=(\underbrace{1,\dots,1}_{m_2\text{
times}}, \underbrace{2,\dots,2}_{m_3\text{ times}},
\underbrace{3,\dots,3}_{m_4\text{ times}}\dots).
\end{equation}
Firstly, observe that all the numbers $(a_s),(b_s),(c_s),(d_s)$ are
uniquely determined by $a_1$ and the collection of increments
\begin{multline}
\label{eq:increment} (a_{s+1}-a_s)_{s=1,\dots,p-1} , (b_1-a_p),
(b_{s+1}-b_s)_{s=1,\dots,q-1},\\  (c_p-b_q),
(c_s-c_{s+1})_{s=1,\dots,p-1},
              (d_q-c_1) , (d_s-d_{s+1})_{s=1,\dots,q-1}, (a_1-d_1),
\end{multline}
where numbers \eqref{eq:increment} are positive integers the sum of
which is equal to $n$.

Every partition of the form depicted on Figure \ref{fig:krzyz} has
free index equal to $(-1)$. There are $n$ choices of $a_1$. There
are $\binom{m_2+m_3+\cdots}{m_2,m_3,\dots}$ different permutations
of the tuple \eqref{eq:gromada}. We need to specify $1\leq p\leq
m_2+m_3+\cdots$. Now we count in how many different ways every of
the numbers \eqref{eq:zlepyciagi} can be written as a sum of two
specified elements of \eqref{eq:increment} minus $1$ except for one
of the numbers \eqref{eq:zlepyciagi} which should be written as a
sum of four specified elements of \eqref{eq:increment} minus $3$;
the position of this exceptional number is specified by the index
$p$. Finally, we have to take into account the symmetry factor
$\frac{1}{4}$ of Figure \ref{fig:krzyz} since for every partition of
this form we can choose sequences $(a_s),(b_s),(c_s),(d_s)$ in four
different ways corresponding to four cyclic rotations of the
template. From Lemma \ref{lem:summertime} it follows that the
contribution of partitions from Figure \ref{fig:krzyz} to the
summand $\Sigma_{\underbrace{1,\dots,1}_{m_2\text{ times}},
\underbrace{2,\dots,2}_{m_3\text{ times}},\dots}$ is equal to
\begin{equation}
\label{eq:wkladkrzyza}(-1) \frac{1}{4} n
\binom{m_2+m_3+\cdots}{m_2,m_3,\dots} \sum_{1\leq p \leq
m_2+m_3+\cdots} \binom{k'_p+2}{3} \prod_{p'\neq p}
\binom{k'_{p'}}{1}.
\end{equation}

By very similar considerations one can easily show that the
contribution of partitions from Figure \ref{fig:gwiazdadawida} to
the summand $\Sigma_{\underbrace{1,\dots,1}_{m_2\text{ times}},
\underbrace{2,\dots,2}_{m_3\text{ times}},\dots}$ is equal to
\begin{multline}
\label{eq:wkladgwiazdydawida}
(-2)\frac{1}{6} n \binom{m_2+m_3+\cdots}{m_2,m_3,\dots} \times \\
\sum_{1\leq p_1<p_2 \leq m_2+m_3+\cdots} \binom{k'_{p_1}+1}{2}
\binom{k'_{p_2}+1}{2} \prod_{p'\neq p_1,p_2} \binom{k'_{p'}}{1}.
\end{multline}

By an application of Leibnitz rule we see that the sum of
\eqref{eq:wkladkrzyza} and \eqref{eq:wkladgwiazdydawida} can be
written in a compact form as
\begin{multline*}
 [x^{-(2m_2+3m_3+\cdots)-2} z_2^{m_2} z_3^{m_3} \cdots] \frac{(-1) n}{24}
\times \\
\shoveright{\frac{d^2}{dx^2} \left( \sum_{s}  (s-1) z_s x^{-s}
\right)^
{m_2+m_3+\cdots}=} \\
[z_2^{m_2} z_3^{m_3} \cdots] \frac{(-1)n(n-1)(n-2)}{24}
\left( \sum_{s}  (s-1) z_s \right)^
{m_2+m_3+\cdots}=\\
\frac{(-1)n(n-1)(n-2)}{24} \binom{m_2+m_3+\cdots}{m_2,m_3,\dots}
1^{m_2} 2^{m_3} 3^{m_4} \cdots
\end{multline*}
which finishes the proof.
\end{proof}

\subsubsection{Second--order expansion of characters}
The following result was conjectured by Biane
\cite{BianeCharacters}.
\begin{theorem}
\label{theo:hahahauuau} For every $n\geq 2$ we have
\begin{multline*}
\Sigma_{n-1}= \RJM_{n}+\\ \!\!\!\!  \sum_{\substack{m_2,m_3,\dots\geq 0 \\
2m_2+3 m_3+\cdots=n-2}} \!\!\!\! \frac{1}{4} \binom{n}{3}
\binom{m_2+m_3+\cdots}{m_2,m_3,\dots} \prod_{s\geq 2} \big( (s-1)
\RJM_s \big)^{m_s} +\\ \text{(terms of degree at most $n-4$)}.
\end{multline*}
\end{theorem}
\begin{proof}
In \eqref{eq:dukakao} we use \eqref{eq:eksplozjamozgu} to express
$\Sigma_{\underbrace{1,\dots,1}_{m_2 \text{
times}},\underbrace{2,\dots,2}_{m_3 \text{ times}},\dots}$ in terms
of $\RJM_2,\RJM_3,\dots$.
\end{proof}

\subsection{Higher order expansion of characters}
In principle, there are no obstacles to repeat the reasoning from
Section \ref{sec:secondorderexp} for partitions with genus $2$ and
in this way obtain higher order asymptotic expansions.

The first step is to enumerate all possible outcomes of the
simplification algorithm. Theorem \ref{theo:gorneszacowanie} shows
that such outcome is a pair partition with at most $9$ lines. The
patience of the most human beings is not sufficient enumerate all
such partitions, nevertheless the use of computer allows us to find
them all \cite{Sniady2004table}. After removing cyclic rotations and
mirror images there are $61$ such partitions---which is still an
accessible number. By reverting all three steps of the
simplification algorithm we could therefore find general patterns of
evercrossing partitions with genus $2$ and therefore express
$\RJM_n$ in terms of elements $\Sigma$.

However, the above method seems to be far more complicated than the
result of Goulden and Rattan \cite{GouldenRattan05} who found an
explicit formula for the general coefficients of Kerov polynomials.

\section{Proofs of technical results}
\label{sec:proofs}

\cytat{
\begin{quotation}
Czucie i wiara silniej m\'owi do mnie ni\.z m\c{e}drca szkie\l{}ko i
oko.
\\
\textsc{Adam Mickiewicz, `Romantyczno\'s\'c'}
\end{quotation}
}

\subsection{Admissible and pushing sequences}
In this article we manipulate with moments of the Jucys--Murphy
element $J$ (or with products of such moments) and we need an
efficient way to enumerate the summands in sums similar to
\eqref{eq:wielkasuma}. Conceptually the simplest way is to enumerate
the summands by sequences $a_1,\dots,a_n\in A$ and which will be
called admissible. However, it turns out that manipulating
admissible sequences is quite complicated. It turns out that the
suitable objects to enumerate summands in \eqref{eq:wielkasuma} are
pushing sequences. We will describe pushing sequences in terms of
admissible sequences.


\subsubsection{Admissible sequences. Map $\fass$}
We recall that $\xtra\notin A$. We say that a sequence
$a_1,\dots,a_n$ is admissible \index{admissible sequence}
\index{sequence!admissible} if it contributes in the sum
\eqref{eq:wielkasuma} or equivalently if $a_1,\dots,a_n\in A$ and
\begin{equation}
[(a_1 \xtra) \cdots (a_n \xtra)](\xtra)=\xtra. \label{eq:admissible}
\end{equation}
The set of admissible sequences is a semigroup if for multiplication
we take the concatenation of sequences:
$$\aaa \bbb=(a_1,\dots,a_n,b_1,\dots,b_o)$$
for admissible sequences $\aaa=(a_1,\dots,a_n)$ and
$\bbb=(b_1,\dots,b_o)$. The unit is the empty sequence.

We consider a map $\fassSym$ from the set of admissible sequences to
the symmetric group $\Sn{A}$ by setting
$$\fassSym:(a_1,\dots,a_n)\mapsto (a_1 \xtra) \cdots (a_n \xtra).$$
We can treat the right--hand side as a partial permutation of $A$
with the support $\{a_1,\dots,a_n\}$; in this way we obtain a map
$\fass$ from the set of admissible sequences to the semigroup
$\PartPerm{A}$ of partial permutations. It is easy to check that
both maps are homomorphisms of semigroups.

\subsubsection{Pushing sequences}
%
We say that no neighbor elements of a sequence $(p_1,\dots,p_n)$ are
equal if $p_l\neq p_{l+1}$ holds for any $1\leq l\leq n$. We should
think that the elements of the sequence are arranged in a circle and
the successor of $p_n$ is $p_1$, therefore in the above definition
the case $l=n$ means that we require $p_n\neq p_1$. In particular,
it cannot happen that $n=1$.

We will say that a sequence $(p_1,\dots,p_n)$ is a pushing sequence
\index{pushing sequence} \index{sequence!pushing} if
$p_1,\dots,p_n\in A\cup\{\xtra\}$ and if $p_1=\xtra$ and if no
neighbor elements of $(p_1,\dots,p_n)$ are equal.

\subsubsection{Commutative diagram} In the following we will construct
functions such that the following diagram commutes. Later on we
shall equip pushing sequences with a multiplicative structure in
such a way that all maps will become homomorphisms of semigroups.
\begin{equation}
\includegraphics{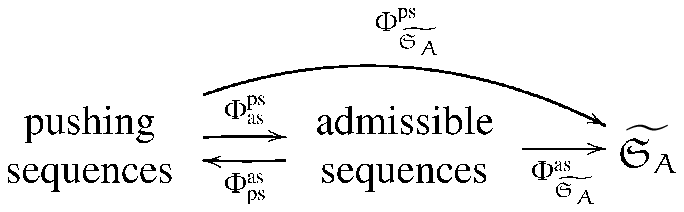}
\label{eq:diagram}
\end{equation}


%

\subsubsection{Map $\fasps$} \label{subsubsec:fasps} Let
$\aaa=(a_1,\dots,a_n)$, $a_1,\dots,a_n\in A$ be an admissible
sequence. We assign to it a sequence of permutations
$\sigma_1,\dots,\sigma_{n+1}\in\Sn{A\cup\{\xtra\}}$ defined by
\begin{equation}
\sigma_l= (a_l \xtra) (a_{l+1} \xtra) \cdots (a_n \xtra)
\label{eq:slynnesigma}
\end{equation}
and a sequence $\ppp=\fasps(\aaa)= (p_1,\dots,p_n)$,
$p_1,\dots,p_n\in A\cup\{\xtra\}$ defined by
\begin{equation}
p_l= \sigma_l^{-1}(\xtra). \label{eq:dlaczegopushing}
\end{equation}
Notice that since $\aaa$ is admissible hence $p_1=\xtra$ and
$p_1=\xtra\neq a_n=p_n$. Furthermore since $\sigma_{l-1}=(a_{l-1}
\xtra) \sigma_l$ hence $p_{l-1}=\sigma_l^{-1}(a_{l-1}) \neq
\sigma_l^{-1}(\xtra)=p_l$ and $\ppp$ is a pushing sequence.

Equation \eqref{eq:dlaczegopushing} should explain the name of
pushing sequences: imagine that $A\cup\{\xtra\}$ is a set of some
items and initially every item $x$ is located within a box which
carries a label `$x$'. In the process of computing the product $(a_1
\xtra) \cdots (a_n \xtra)$ (from the right to the left) we obtain a
sequence of partial products $\sigma_n,\sigma_{n-1},\dots,\sigma_1$:
permutation $\sigma_{l-1}$ is obtained from $\sigma_l$ by
transposition of the item which is in the box labeled by `$\xtra$'
and the item (which turns out to be the item $p_{l-1}$) from the box
labeled by `$a_{l-1}$'. We can think therefore that the item
$p_{l-1}$ is pushing out the contents of the box labeled by
`$\xtra$'.


\subsubsection{Map $\fpsas$} \label{subsubsec:fpsas}

Let a pushing sequence $\ppp=(p_1,\dots,p_n)$ be given. By backward
induction we will assign to it a unique pair of sequences
$(\sigma_1,\dots,\sigma_{n+1})$, $\fpsas(\ppp)= (a_1,\dots,a_n)$,
where $\sigma_1,\dots,\sigma_{n+1}\in\Sn{A\cup\{\xtra\}}$ and
$(a_1,\dots,a_n)$ is an admissible sequence such that
\eqref{eq:slynnesigma} and \eqref{eq:dlaczegopushing} are fulfilled.

We start the construction by setting $\sigma_{n+1}=e$ since it is
the only choice if we want \eqref{eq:slynnesigma} to be fulfilled
for $l=n+1$. If $\sigma_l$ is already defined we set
$a_{l-1}=\sigma_{l} (p_{l-1})$ and $\sigma_{l-1}=(a_{l-1} \xtra)
\sigma_l$. Again, it is the only choice that we have since equations
 \eqref{eq:slynnesigma} and \eqref{eq:dlaczegopushing} imply
\begin{equation}
\label{eq:formulanaadmissible} a_{l-1}=(a_{l-1}\xtra)(\xtra)=
\sigma_{l} \sigma_{l-1}^{-1}(\xtra)=\sigma_{l}(p_{l-1}).
\end{equation}

Since $\ppp$ is pushing hence $\sigma_l^{-1}(a_{l-1})=p_{l-1} \neq
p_l=\sigma_l^{-1}(\xtra)$ therefore $a_1,\dots,a_n\neq \xtra$;
furthermore $[(a_1 \xtra) \cdots (a_n
\xtra)]^{-1}(\xtra)=\sigma_1^{-1}(\xtra)=p_{1}=\xtra$ hence $\aaa$
is admissible.

Equations \eqref{eq:slynnesigma} and \eqref{eq:dlaczegopushing} are
indeed fulfilled by construction, therefore
$\fasps(\fpsas(\ppp))=\ppp$ for every pushing sequence $\ppp$.
Furthermore the uniqueness of the pair $(a_1,\dots,a_n)$,
$(\sigma_1,\dots,\sigma_{n+1})$ implies that there exists at most
one preimage of $\ppp$ with respect to $\fasps$. Therefore we proved
the following result:
\begin{proposition}
Maps $\fpsas$ and $\fasps$ are inverses of each other.
\end{proposition}




%
%

\subsubsection{Map $\fpss$} The bijections between pushing and
admissible sequences allow us to define the map $\fpss$ in a unique
way which makes the diagram \eqref{eq:diagram} commute.

Let $\ppp$ be a pushing sequence and $\aaa$ be the corresponding
admissible sequence. Observe that if $\sigma_1,\dots,\sigma_{n+1}$
is the sequence of permutations which was constructed in Sections
\ref{subsubsec:fasps} and \ref{subsubsec:fpsas} then
\begin{equation}
\label{eq:formulanafpss} \fpssSym(\ppp)=\fassSym(\aaa)=\sigma_1.
\end{equation}
In order to compute $\fpss(\ppp)$ we need to specify the support of
this permutation. It is the support of $\fass(\aaa)$, namely
$\{a_1,\dots,a_n\}=\{p_1,\dots,p_n\}\setminus\{\xtra\}$; the
inclusions which imply the latter equality follow directly from the
construction of the maps $\fasps$ and $\fpsas$. In Section
\ref{subsec:explicit} we shall compute $\fpss$ and $\fpsas$ more
explicitly.

\subsubsection{Multiplication of pushing sequences}
\label{subsec:multiplicationpushing} The bijections between pushing
and admissible sequences allow us to define multiplication of
pushing sequences in the unique way which makes arrows of
\eqref{eq:diagram} to be homomorphisms of semigroups in such a way
that the diagram commutes.

The exact form of this multiplication is given by the following
proposition
\begin{proposition}
\label{prop:mnozeniepushingsequences} Let $\ppp=(p_1,\dots,p_k)$ and
$\qqq=(q_1,\dots,q_l)$ be pushing sequences and let
$\pi=\fpssSym(q_1,\dots,q_l)$. Then
$$\ppp \qqq= \big(\pi^{-1}(p_1),\dots,\pi^{-1}(p_k),q_1,\dots,q_l \big). $$
\end{proposition}
\begin{proof}
Let $\ppp=\fasps(\aaa)$ and let $\sigma^a_1,\dots,\sigma^a_k$ be the
sequence of permutations constructed in Section
\ref{subsubsec:fasps} for $\aaa$, i.e.~$\sigma^a_m=(a_m
\xtra)\cdots(a_k\xtra)$ for $1\leq m\leq k$. Let $\bbb=\fasps(\qqq)$
and let $\sigma^b_1,\dots,\sigma^b_l$ be the analogous sequence of
permutations obtained for $\bbb$.

Let $\sigma_1,\dots,\sigma_{k+l}$ be the sequence of permutations
obtained for the concatenated sequence $\aaa \bbb$. Observe that
\begin{multline*}\sigma_m=\begin{cases}
(b_{m-k} \xtra) \cdots (b_l \xtra) & \text{for } m\geq k+1 \\
(a_{m} \xtra) \cdots (a_k \xtra) (b_1 \xtra) \cdots (b_l \xtra) &
\text{for } m\leq k \end{cases}  \\ =
\begin{cases} \sigma^b_{m-k}  & \text{for } m\geq k+1 \\
\sigma^a_{m} \pi & \text{for } m\leq k,
\end{cases}
\end{multline*}
where $\pi=(b_1 \xtra) \cdots (b_l
\xtra)=\fassSym(\bbb)=\fpssSym(\qqq)$.

It follows that $\rrr=(r_1,\dots,r_{n+l})= \ppp \qqq=\fasps(\aaa
\bbb)$ is given by
$$r_m=\sigma_m^{-1}(\xtra)=\begin{cases} q_{m-n} & \text{for } m\geq n+1 \\
\pi^{-1}(p_{m}) & \text{for } m\leq n
\end{cases} $$
which finishes the proof.
\end{proof}

\subsubsection{Invariance of the conjugacy class assigned to a pushing
sequence}
For the purpose of this section let us forget about the set $A$ and
the distinguished element $\xtra$. Let a sequence
$\ppp=(p_1,\dots,p_n)$ be given such that no neighbor elements are
equal: $p_l\neq p_{l+1}$ for all $1\leq l\leq n$. Now we can define
the distinguished element $\xtra$ by setting $\xtra=p_1$ and set $A$
by $A=\{p_1,\dots,p_n\}\setminus \{\xtra\}$.

Let $f$ be a one to one function defined on the set
$\{p_1,\dots,p_n\}$. We set $\ppp'=\big( f(p_1),\dots,f(p_n) \big)$,
$\xtra'=f(\xtra)$, $A'=f(A)$.

Pushing sequence $\ppp$ defines a partial permutation
$\Phi^{\ps}_{\PartPerm{A}}(\ppp)\in \PartPerm{A}$ and the pushing
sequence $\ppp'$ defines a partial permutation
$\Phi^{\ps}_{\PartPerm{A'}}(\ppp')\in \PartPerm{A'}$.

\begin{proposition}
\label{prop:zmianaoznaczen} We use the above notations. Then
$$\fpssSym(\ppp')=f\circ \fpssSym(\ppp) \circ f^{-1}$$
and the support of $\Phi^{\ps}_{\PartPerm{A'}}(\ppp')$ is equal to
the image of the support of $\Phi^{\ps}_{\PartPerm{A}}(\ppp)$ under
the map $f$.
%
\end{proposition}
\begin{proof}
Let $a_1,\dots,a_n$ and $\sigma_1,\dots,\sigma_{n+1}$ be the
sequences we considered in Sections \ref{subsubsec:fasps} and
\ref{subsubsec:fpsas} and let $a'_1,\dots,a_n'$ and
$\sigma_1',\dots,\sigma_{n+1}'$ be their analogues obtained for the
pushing sequence $\ppp'$.

We leave to the reader to show by backward induction that
$a_l'=f(a_l)$ and $\sigma_l'= f\circ \sigma_l \circ f^{-1}$ for
every $1\leq l\leq n$. From \eqref{eq:formulanafpss} (and its
analogue for $\ppp'$) it follows that
$$\fpssSym(\ppp')=f\circ \fpssSym(\ppp) \circ f^{-1}.$$

Now it is enough to observe that the support of
$\Phi^{\ps}_{\PartPerm{A'}}(\ppp')$ is equal to $A'$ and the support
of $\Phi^{\ps}_{\PartPerm{A}}(\ppp)$ is equal to $A$.
\end{proof}

\subsection{Pushing partitions}
In this article we are concerned mostly with conjugacy classes of
permutations and not permutations themselves. For this reason we
need a method to enumerate in the sums similar to
\eqref{eq:wielkasuma} not the individual summands, but whole classes
of summands (summands within each class will be conjugate). Pushing
partitions will turn out to be the suitable objects to enumerate
such classes.

We say that a partition $\pi$ of a finite ordered set $X$ (by
changing the labels we can assume that $X=\{1,\dots,n\}$) is a
pushing partition if no neighbors $l$ and $l+1$ are connected by
$\pi$ for $1\leq l\leq n$. It is useful to arrange numbers
$1,\dots,n$ in a circle, therefore the case $l=n$ should be
understood that neighbors $n$ and $1$ are not connected by $\pi$. In
particular, it cannot happen that $n=1$.

We say that $\ppp\sim\pi$ if $\pi$ is a partition of $\{1,\dots,n\}$
for some $n\in\N$ and $\ppp=(p_1,\dots,p_n)$ is a sequence of length
$n$ such that for every $1\leq l,m\leq n$ numbers $l$ and $m$ are
connected by $\pi$ iff $p_l=p_m$. For every pushing sequence $\ppp$
there exists a unique pushing partition $\pi$ such that
$\ppp\sim\pi$.


\subsubsection{Explicit form of the maps $\fpss$ and $\fpsas$}
\label{subsec:explicit} Let a pushing partition
$\pi=\{\pi_1,\dots,\pi_r\}$ of the set $\{1,2,\dots,n\}$ be given.
To every block $\pi_s$ of the partition we assign a label $\min
\pi_s$, i.e.~the smallest element of $\pi_s$. To every element of
$\{1,2,\dots,n\}$ we assign the label of the block it belongs to. In
this way we constructed a sequence $\ppp=(p_1,\dots,p_n)$ such that
$\ppp \sim \pi$. We set $\xtra=p_1=1$ and $A$ to be the set of the
other labels, i.e.~$A=\{p_1,\dots,p_n\}\setminus\{\xtra\}=\{\min
\pi_s: 1\leq s\leq r\}\setminus \{1\}$. We decorate on the graph of
the fat partition $\pi_{\fat}$ all elements of the set $A$ (cf
Figure \ref{fig:tlustapartycja2}).

In the following we shall compute $\fpss(\ppp)$ for $\ppp$ defined
as above; please note that Proposition \ref{prop:zmianaoznaczen}
allows us to rename the labels and to find $\fpss(\ppp)$ for a
general $\ppp$.

As we already mentioned in Section
\ref{subsec:pierwszawzmiankaofppsq}, we can view the fat partition
$\pi_{\fat}$ as a bijection
$\pi_{\fat}:\{1',2',\dots,n'\}\rightarrow\{1,2,\dots,n\}$. We also
consider a bijection
$c:\{1,2,\dots,n\}\rightarrow\{1',2',\dots,n'\}$ given by
$\dots,3\mapsto 2', 2\mapsto 1', 1\mapsto n', n\mapsto
(n-1)',\dots$. We consider the cycle decomposition of the
permutation $\pi_{\fat}\circ c\in\Sn{\{1,2,\dots,n\}}$ and remove
from this decomposition all elements except the elements of $A$. In
this way we constructed a permutation $(\pi_{\fat}\circ
c)|_A\in\Sn{A}$, called restriction of $\pi_{\fat}\circ c$ to the
set $A$. This operation has a nice interpretation: on the graphical
representation of $\pi_{\fat}$ and $c$ we travel along the cycles by
following the arrows and we write down only decorated vertices. We
equip $(\pi_{\fat}\circ c)|_A$ with a structure of a partial
permutation by setting its support to be $A$.

For example, for $\pi=\big\{ \{1,3\},\{2,5,7\}, \{4\}, \{6\} \big\}$
depicted on Figure \ref{fig:przecinajacapartycja} we have
$\ppp=(1,2,1,4,2,6,2)$, $A=\{2,4,6\}$ and the decorations of
$\pi_{\fat}$ are depicted on Figures \ref{fig:tlustapartycja2},
\ref{fig:tlustapartycja3} and \ref{fig:tlustapartycja4}. The
composition $\pi_{\fat}\circ c$ has a cycle decomposition
$(1,\mathbf{2},3,5,\mathbf{4})(\mathbf{6},7)$.
It follows that $(\pi_{\fat}\circ
c)|_A=(\mathbf{2},\mathbf{4})(\mathbf{6})$.

\begin{theorem}
\label{theo:ilewynosifpss} Let $\pi$ and $\ppp$ be as above. Then
$$ \fpss(\ppp)=(\pi_{\fat} \circ c)|_A.$$
\end{theorem}
\begin{proof}
The equality of permutations follows easily from Lemma
\ref{lem:lematoliczniufpss} below by setting $l=1$ and the equality
of supports is trivial.
\end{proof}
\begin{lemma}
\label{lem:lematoliczniufpss} Let $\pi$ and $\ppp$ be as above, we
set $p_{n+1}:=p_1=1=\xtra$ and let $\sigma_1,\dots,\sigma_{n+1}$ be
as in Sections \ref{subsubsec:fasps} and \ref{subsubsec:fpsas}. Let
$x\in A\cup\{\xtra\}$ and let $1\leq l\leq n+1$.

If $x=p_l$ then $\sigma_l(x)=\xtra$.

Suppose that $x\neq p_l$; let $m$ be the smallest index such that
$l\leq m\leq n+1$ and $x=p_m$. If no such index exist then
$\sigma_l(x)=x$. Otherwise we start in the vertex $m$ a walk on the
graphical representation of $\pi_{\fat}$ and $c$ by following the
arrows until we enter some decorated vertex $y$ (after having made
at least one step, i.e.~starting in a decorated vertex do not count
as entering it). The above walk always stops after a finite number
of steps and $\sigma_l(x)=y$.
\end{lemma}
\begin{proof}
The shall prove the lemma by backward induction. In the case $l=n+1$
it holds trivially.

Suppose that the statement of lemma is true for all $l'>l$. The
first part, $\sigma_l(p_l)=\xtra$ follows from the very construction
of $\sigma_l$ in \eqref{eq:dlaczegopushing}.

Observe that for every $r$ we have $\sigma_r=(a_r\xtra)
\sigma_{r+1}$ hence $\sigma_r(x)=\sigma_{r+1}(x)$ if $x\notin\{
\sigma_{r}^{-1}(\xtra),
\sigma_{r+1}^{-1}(\xtra)\}=\{p_{r},p_{r+1}\}$. This has twofold
implications.

Firstly, if $x\notin\{p_l,p_{l+1},\dots,p_n,p_{n+1}\}$ then
$\sigma_l(x)=\sigma_{l+1}(x)=\cdots=\sigma_{n+1}(x)=x$, which
finishes the proof.

Secondly, when the last case in the lemma statement holds then
$$\sigma_l(x)=\sigma_{l+1}(x)=\cdots=\sigma_{m-1} (x)=
  \big((a_{m-1} \xtra)\circ \sigma_m\big)(x)=a_{m-1},$$
where we used the fact that $p_m=x$ is equivalent to
$\sigma_m(x)=\xtra$. Equation \eqref{eq:formulanaadmissible} implies
that $\sigma_l(x)=\sigma_{m}(p_{m-1})$. We denote by $\pi_h$ the
block of $\pi=\{\pi_1,\dots,\pi_r\}$ such that the element $m-1$
belongs to it. Let $r$ be the smallest index such that $p_{m-1}=p_r$
and $m\leq r\leq n+1$. If no such index exists then
$\sigma_l(x)=\sigma_{m}(p_{m-1})=p_{m-1}$. Otherwise, if such an
index exists then
$$\sigma_l(x)=\sigma_{m}(p_{m-1})=\sigma_{m+1}(p_{m-1})=\cdots=
\sigma_{r}(p_{m-1})=\sigma_{r}(p_r).$$

We compare the value of $\sigma_l(x)$ computed above with the answer
given by the graphical algorithm: we start our trip in vertex $m$
and then go to vertex $(m-1)'$. If the index $r$ considered above
does not exist then $m-1$ is the biggest element of $\pi_h$ and
therefore $(m-1)'$ is connected by $\pi_{\fat}$ with the smallest
element of $\pi_h$ which is therefore decorated (the exceptional
case $1\in \pi_h$ is not possible since this would mean that
$p_{m-1}=p_1=\xtra$ and in this case $p_{m-1}=p_{n+1}$ and the index
$r$ would exist) and the algorithm terminates giving the answer
$\min \pi_h=p_{m-1}$ which coincides with the correct value.
Otherwise, if the index $r$ exists then  $(m-1)'$ is connected by
$\pi_{\fat}$ with the vertex $r$, $l<m+1\leq r\leq n+1$ (the vertex
$n+1$ should be understood as the vertex $1$). From the inductive
hypothesis it follows that the algorithm will return the answer
$\sigma_{r}(p_r)$ which is again correct.
\end{proof}

\begin{corollary}
\label{cor:formulanaadmissible} Let $\pi$ and $\ppp$ be as above,
let $\aaa=\fpsas(\ppp)$ and let $1\leq l\leq n$. We start in the
vertex $l$ a walk on the on the graphical representation of
$\pi_{\fat}$ and $c$ by following the arrows until we enter some
decorated vertex $y$ (after having made at least one step).

Then
$$a_l=y.$$
\end{corollary}
\begin{proof}
It is a simple application of \eqref{eq:formulanaadmissible} and
Lemma \ref{lem:lematoliczniufpss}.
\end{proof}

\subsubsection{Maps $\fppps$ and $\Sigma=\fpps$} Let us fix the finite set
$A$ and the extra element $\xtra$. We consider an algebra of pushing
sequences in which the multiplication is as in Section
\ref{subsec:multiplicationpushing} and the addition is understood as
an addition of formal sums. The map $\fpss$ extends naturally to an
algebra homomorphism from the algebra of pushing sequences to the
algebra $\C(\PartPerm{A})$ of partial permutations.

We define a map from pushing partitions to the algebra of pushing
sequences by
\begin{equation}
\label{eq:definicjafppps} \fppps(\pi)=\sum_{\ppp\sim\pi} \ppp.
\end{equation}

We define the map $\fpps$ in the unique way which makes the diagram
below commute. This map is explicitly given by Claim
\ref{claim:definicjasigma} which will be proved below.
\begin{equation}
\label{eq:diagram2}
\includegraphics{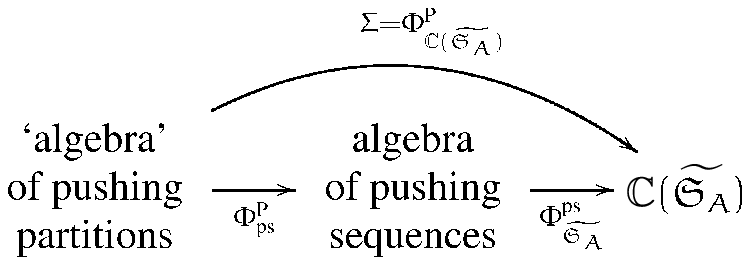}
\end{equation}

\begin{theorem}
\label{theo:stwierdzeniepierwszejestprawda} Claim
\ref{claim:definicjasigma} is true. In other words: let $\pi$ be a
partition of the set $\{1,2,\dots,n\}$ and let numbers
$k_1,\dots,k_t$ be given by the above construction. Then
$$\Sigma_\pi:=\fppsinf(\pi)=\Sigma_{k_1,\dots,k_t},$$
where $\Sigma_{k_1,\dots,k_t}\in\C(\PartPerm{\infty})$ on the
right--hand side should be understood as in Section
\ref{subsec:definicjasigma}.
\end{theorem}
\begin{proof}
Equation \eqref{eq:definicjafppps} and requirement that
\eqref{eq:diagram2} commutes imply that
$$\Sigma(\pi)=\sum_{\ppp\sim\pi} \fpss(\ppp).$$
By an appropriate renaming of labels Theorem
\ref{theo:ilewynosifpss} allows us to compute each summand
$\fpss(\ppp)$. We decompose $\fpss(\ppp)=(\pi\circ c)|_A
\in\Sn{\{p_1,p_2,\dots,p_n\} \setminus\{p_1\} }$ into a product of
disjoint cycles and we denote by $k_1,\dots,k_t$ the lengths of
these cycles. It follows that
$$\Sigma_{\pi}=\Sigma_{k_1,\dots,k_t},$$
where the right--hand side was defined in \eqref{eq:definicjasigma}.

From the proof of Theorem \ref{theo:ilewynosifpss} it follows that
for each $1\leq s\leq t$ the number $k_s$ is equal to the number of
decorated elements in the corresponding cycle of the permutation
$\pi_{\fat}\circ c$. One can easily see that the original definition
of the decorated elements is equivalent to the following one:
element $m\in\{1,2,\dots,n\}$ is decorated if and only if going
counterclockwise from $m$ to $(\pi_{\fat}\circ c)^{-1}(m)$ one does
not cross the line between the marked starting point between $1$ and
$1'$ and the central disc (see Figure \ref{fig:tlustapartycja2}).
Therefore $k_s$ is equal to the number of elements in the
corresponding cycle of $\pi_{\fat}\circ c$ minus the number of times
this cycle going clockwise crosses the line between the starting
point and the central disc. It follows that the numbers
$k_1,\dots,k_t$ defined above coincide with \eqref{eq:formulanak}
which finishes the proof.
\end{proof}

\subsubsection{Multiplication of pushing partitions} In Section
\ref{sec:mnozeniepartycji} we defined the multiplication of
partitions. Diagram \eqref{eq:diagram2} commutes and we shall prove
that $\Sigma=\fpps$ is an algebra homomorphism (we already know that
$\fpss$ is an algebra homomorphism; we will not need it but one can
show that in a certain sense $\fppps$ is also a homomorphism).


\begin{lemma}
\label{lem:mnozeniedwoch} Let $1\leq k< l\leq m$, let $\pi^1$ be a
pushing partition of the set
$\rho_1=\{l+1,l+2,\dots,m,1,2,\dots,k\}$, let $\pi^2$ be a pushing
partition of the set $\rho_2=\{k+1,k+2,\dots,l\}$. We define a
non--crossing partition $\rho=\{\rho_1,\rho_2\}$.

Then $$\fpps(\pi^1\cdot \pi^2)=\fpps(\pi^1) \fpps(\pi^2),$$ where
$\pi^1 \cdot \pi^2$ denotes the $\rho$--ordered product of
partitions $\pi^1$ and $\pi^2$ as defined in Section
\ref{sec:mnozeniepartycji} and the product on the right--hand side
is the usual product of the elements of $\C(\PartPerm{A})$.
\end{lemma}
\begin{proof}
Let us consider a simpler case $\rho'_1=\{1,2,\dots,k+m-l\}$,
$\rho'_2=\{k+m-l+1,k+m-l+2,\dots,m\}$. From Proposition
\ref{prop:mnozeniepushingsequences} it follows that $\fppps(\pi^1)
\fppps(\pi^2)$ consists of all pushing sequences $(p_1,\dots,p_m)$
such that $p_1=p_{k+m-l+1}=\xtra$ and $(p_1,\dots,p_{k+m-l})\sim
\pi^1$, $(p_{k+m-l+1},\dots,p_m)\sim \pi^2$. On the other hand
the only non--trivial block of $\rho'_{\comp^{-1}}$ is
$\{1,k+m-l+1\}$ and therefore
%
%
$$\fppps(\pi^1) \fppps(\pi^2) = \fppps( \pi^1
\cdot \pi^2).$$ To both sides of the equation we apply the map
$\fpss$ and use the fact that it is a homomorphism and the diagram
\eqref{eq:diagram} commutes.

The general case follows from the observation that the partition
$\rho$ can be obtained from $\rho'$ by a cyclic rotation. One can
easily see that a product of rotated partitions is a rotation of the
product of the original partitions. Finally, we apply Proposition
\ref{prop:rotacjajestok}.
\end{proof}

\begin{lemma}
\label{lem:associative} Multiplication of partitions, as defined in
Section \ref{sec:mnozeniepartycji}, is associative. To be precise:
let $\tau=\{\tau_1,\dots,\tau_r\}$ be a non--crossing partition of
the set $\{1,2,\dots,n\}$. For each $s$ let
$\tilde{\tau}^s=\{\tau_{s,1},\dots,\tau_{s,l_s} \}$ be a
non--crossing partition of the set $\tau_s$. We denote by
$\hat{\tau}=\tilde{\tau}^1\cup\cdots\cup\tilde{\tau}^r = \{
\tau_{s,k}: 1\leq s\leq r, 1\leq k\leq l_s  \}$ the non--crossing
partition of the set $\{1,2,\dots,n\}$. Let furthermore $\pi^{s,k}$
be a pushing partition of the set $\tau_{s,k}$ for each choice of
$1\leq s\leq r$ and $1\leq k\leq l_s$.

Then
\begin{equation}
 \prod_s \Big(  \prod_k   \pi^{s,k}  \Big)= \prod_{s,k}   \pi^{s,k},
\label{eq:nudynapudy}
\end{equation}
where the left--hand side is a $\tau$--ordered product of
$\tilde{\tau}^s$--ordered products and the right--hand side is a
$\hat{\tau}$--ordered product.
\end{lemma}
\begin{proof}
Probably the best way to prove this lemma is to use the graphical
description of the partition multiplication from Section
\ref{subsec:geometricmultiplication} since---speaking
informally---cutting holes, gluing discs and merging discs are all
associative operations. We provide a more detailed proof below.

Observe that every partition $\sigma$ which contributes to the
left--hand side has a property that for every $a,b\in \tau_{s,k}$
elements $a$ and $b$ are connected by $\sigma$ if and only if they
are connected by $\pi^{s,k}$. Furthermore since every partition
$\sigma^s$ which contributes in $\prod_{k} \pi^{s,k}$ must fulfill
$\sigma^s \geq (\tilde{\tau}^s)_{\comp^{-1}}$ therefore $\sigma\geq
\tau_{\comp^{-1}}\vee\Big(\bigcup_s
(\tilde{\tau}^s)_{\comp^{-1}}\Big)=\hat{\tau}_{\comp^{-1}}$ (for the
last equality see Lemma \ref{lem:babel} below). It follows that
every partition which appears on the left--hand side of
\eqref{eq:nudynapudy} appears also on the right--hand side.

To show the opposite inclusion, for a summand $\sigma$ which appears
on the right--hand side of \eqref{eq:nudynapudy} we set $\sigma^s\in
\PP(\tau_s)$ to be the partition which connects $a,b\in \tau_s$ if
and only if $a,b$ are connected by $\sigma$. Lemma \ref{lem:babel}
implies that $\sigma^s\geq (\tilde{\tau}^s)_{\comp^{-1}}$ hence
$\sigma^s$ is one of the summands which appear in the product
$\prod_k \pi^{s,k}$. By using a similar argument we see that
$\sigma$ is one of the summands which contribute to the product
$\prod_s \sigma^s$ and therefore contributes to the left--hand side
of \eqref{eq:nudynapudy}.
\end{proof}

\begin{lemma}
\label{lem:babel} We keep the notation from Lemma
\ref{lem:associative}. Then
\begin{equation}
\tau_{\comp^{-1}}\vee\Big(\bigcup_s
(\tilde{\tau}^s)_{\comp^{-1}}\Big)=\hat{\tau}_{\comp^{-1}}.
\label{eq:widzeciemnosc}
\end{equation}
\end{lemma}
\begin{proof}
Similarly as in Section \ref{subsec:geometricmultiplication} we
consider a sphere with a collection of holes, each corresponding to
one of the blocks of $\tau$. In the second step we replace each
initial hole $\tau_s$ by a collection of holes $\tilde{\tau}^s$: in
this way we obtained a large sphere with a collection of holes
$\tau_{s,k}$.  Alternatively we can treat it as one big hole to
which was glued a collection of discs; every of these discs
corresponds to one of the blocks of $\hat{\tau}_{\comp^{-1}}$ which
corresponds to the right--hand side of \eqref{eq:widzeciemnosc}.

On the other hand, we can treat the sphere with the first collection
of holes as a sphere with a single hole glued with a collection of
discs corresponding to $\tau_{\comp^{-1}}$. When we replace the hole
$\tau_s$ by $\tilde{\tau}^s$ it corresponds to gluing another
collection of discs $(\tilde{\tau}^s)_{\comp^{-1}}$. Therefore the
whole collection of discs corresponds to the left--hand side of
\eqref{eq:widzeciemnosc}.
\end{proof}

\begin{theorem}
\label{theo:sigmajesthomomorfizmem} Claim \ref{claim:mnozenie} is
true. In other words: let $\rho=\{\rho_1,\dots,\rho_r\}$ be a
non--crossing partition of the set $\{1,2,\dots,n\}$ and for every
$1\leq s\leq r$ let $\pi^s$ be a partition of the set $\rho_s$. Then
$$\Sigma\Big( \prod_{s} \pi^s \Big)=\prod_{1\leq s\leq r} \Sigma_{\pi^s},$$
where the multiplication on the left hand side should be understood
as the $\rho$--ordered product of partitions and on the right hand
side it should be understood as the usual product of commuting
elements in $\C(\PartPerm{\infty})$.
\end{theorem}
\begin{proof}
We shall use the induction with respect to $r$, the number of blocks
of the partition $\rho$. The case $r=1$ is of course trivial.

Observe that since $\rho$ is non--crossing therefore at least one
block $\rho_h$ has a form $\rho_h=\{a,a+1,a+2,\dots,b\}$ (for
example, it is one of the blocks for which the expression $(\max
\rho_h- \min \rho_h)$ takes the minimal value).
Lemma \ref{lem:associative} implies that
$$\prod_s \pi^s = \pi^h \cdot \big( \prod_{s\neq h} \pi^s \big). $$
We apply to both sides the map $\fpps=\Sigma$ and apply Lemma
\ref{lem:mnozeniedwoch} to the right hand side. The inductive
hypothesis can be applied to the second factor on the right hand
side which finishes the proof.
\end{proof}

\subsubsection{Jucys--Murphy element}
\begin{theorem}
\label{theo:zgodnoscpartycjizjurcysiem} Claim
\ref{claim:dobryjurcysmurphy} is true. In other words, let $X$ be a
finite ordered set, let $\rho$ be a non--crossing partition of a
finite ordered set and $n\geq 1$ be an integer. Then \begin{align*}
 \Sigma(\Mpp_{X}) &=\MJM_{|X|}, \\
 \Sigma(\Mpp_\rho)&=\MJM_\rho, \\
 \Sigma(\Rpp_n)&=\RJM_n, \end{align*} where
$\Sigma$ is the map considered in Section
\ref{subsec:pierwszawzmiankaofppsq}.
\end{theorem}
\begin{proof}
Observe that every pushing sequence of length $n$ appears in
$\fppps(\Mpp_{\{1,2,\dots,n\}})$ exactly once. The bijection between
pushing and admissible sequences together with the definition
\eqref{eq:wielkasuma} finish the proof.
\end{proof}


\section{Final remarks}

\label{sec:randommatrices}

\subsection{Connection with the work of Biane}
\label{subsec:bianeconnection} The following connection with article
\cite{Biane1998} was pointed out to me by Philippe Biane.


Biane shows \cite{Biane1998} that the matrix $\Gamma$ defined in
\eqref{eq:macierzjucysia} fulfills
\begin{multline*}
\MJM_n=\frac{1}{q+1} \Tr \Gamma^n=\\
\frac{1}{q+1} \sum_{p_1,\dots,p_{n} \in\{1,\dots,q+1\}} \Gamma_{p_1
p_2} \Gamma_{p_2 p_3} \cdots \Gamma_{p_{n-1} p_n} \Gamma_{p_n p_1},
\end{multline*}
where the sum runs only over tuples $(p_1,\dots,p_n)$ such that no
neighbor elements are equal. Under identification $\xtra=q+1$
summands in Biane's formula coincide with pushing sequences; the
only modification is that one does not assume that $p_1=\xtra$.

This connection goes much deeper: in Section 4.3 of \cite{Biane1998}
one considers pairs $(J,\pi)$ where
$J=\{r_1,\dots,r_k\}\subseteq\{1,2,\dots,n\}$ is the set of indices
such that $p_{r_1}=\cdots=p_{r_k}=\xtra$ and $\pi$ is a partition of
the set $\{1,2,\dots,n\}\setminus J$ which tells us which of the
elements in the tuple $(p_1,\dots,p_n)$ are equal. One can easily
see that the union $\{J\} \cup \pi$ is a pushing partition. We leave
it as an exercise to the reader to verify that the conjugacy class
assigned to a pair $(J,\pi)$ by Biane coincides with the conjugacy
class which appears in $\fppsSym(\{J\} \cup \pi)$. Methods presented
in this article can be used therefore to simplify some of the
arguments in the papers \cite{Biane1998,Biane2001approximate}.

%
%
\subsection{Ramified coverings and pushing partitions}
\label{subsec:okounkov}
%

\subsubsection{Ramified covering and an admissible sequence}
\begin{figure}[bt]
\includegraphics{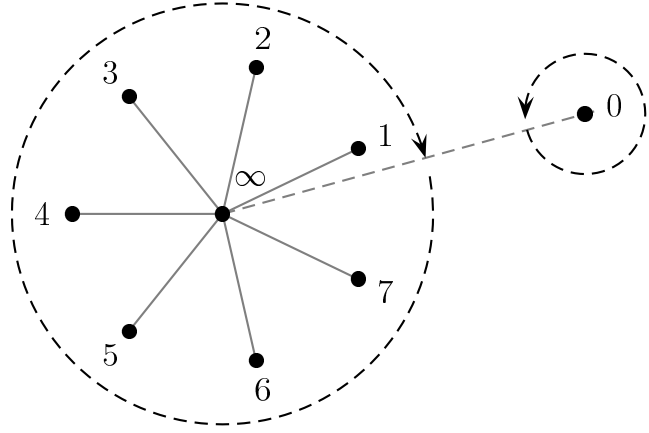}
\caption{Two homotopic loops on a sphere.} \label{fig:sfera}
\end{figure}

Suppose that an admissible sequence
$a_1,\dots,a_n\in\{1,2,\dots,q\}$ is given. Following the idea of
Okounkov \cite{Okounkov2000randompermutations} we consider a
two--dimensional sphere $\Sdwa$ with distinguished points
$0,1,2,\dots,n,\infty$. Figure \ref{fig:sfera} depicts a small
region on this sphere (unfortunately, some conventions used in the
pictures in this article do not coincide with the ones of Okounkov).

Let $\Ss$ be an oriented surface. We shall consider a covering
$\Ss\rightarrow\Sdwa$ with simple ramifications over points
$1,2,\dots,n$ and an unspecified ramification over $0$. We consider
$n+1$ cuts from points $0,1,2,\dots,n$ to $\infty$. We choose $\Ss$
in such a way that in the fiber over any point not lying on a cut
there are $q+1$ sheets marked $1,2,\dots,q,\xtra$. Sheet $\xtra$
will be called special sheet. We also require that the monodromy
around any point $k\in\{1,\dots,n\}$ be a transposition of the sheet
$a_k$ and the sheet $\xtra$. Since a loop around $0$ is homotopic to
a loop around points $1,2,\dots,n$ hence the monodromy around $0$ is
a product of monodromies around $1,\dots,n$ and hence is equal to
$(a_1 \xtra)\cdots(a_n \xtra)$ (cf Figure \ref{fig:sfera}). The
difference with the situation from the paper
\cite{Okounkov2000randompermutations} is that we do not assume that
$(a_1 \xtra)\cdots (a_n \xtra)=e$ and for this reason we might have
a non--trivial ramification in $0$ and therefore we need the extra
cut between $0$ and $\infty$, not present in
\cite{Okounkov2000randompermutations}.

\subsubsection{Shape of the special sheet}
Since the sequence $a_1,\dots,a_n$ is admissible, hence the
monodromy around $0$ preserves the special sheet. It follows that
when we consider only cuts from points $1,2,\dots,n$ to $\infty$
then $\Ss$ splits into the special sheet and the union of
non--special sheets (which are, possibly, glued together along the
cut between $0$ and $\infty$). Let us have a look on the shape of
the special sheet. Along the cuts from $1,2,\dots,n$ to $\infty$ it
is glued with non--special sheets and we are not interested how does
this gluing look like. Much more interesting is the vertex $\infty$
since some points on the boundary of the special sheet might be
glued there together.

\begin{figure}[bt]
\includegraphics{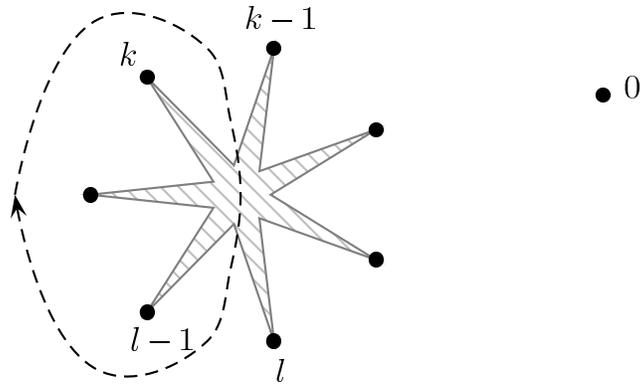}
\caption{A loop on a sphere.} \label{fig:ugaugauga}
\end{figure}

To check if two vertices: the one between $k-1$ and $k$ and the one
between $l-1$ and $l$ are glued together or not, we consider a loop
as on Figure \ref{fig:ugaugauga}. To indicate better the shape of
the special sheet we inflated slightly the cuts. Our question is
equivalent to the following one: is it true that the monodromy along
this loop preserves the special sheet. This holds if and only if
$$\big((a_{k} \xtra)(a_{k+1} \xtra)\cdots (a_{l-1} \xtra)\big)^{-1}(\xtra)=\xtra$$
which is equivalent to
$$p_k=\big((a_{k} \xtra)(a_{k+1} \xtra) \cdots (a_{n} \xtra)\big)^{-1}(\xtra)=
                      \big( (a_{l} \xtra)  \cdots (a_{n} \xtra) \big)^{-1}(\xtra)=p_l,$$
where $p_1,\dots,p_n$ is the corresponding pushing sequence.

\begin{figure}[bt]
\includegraphics{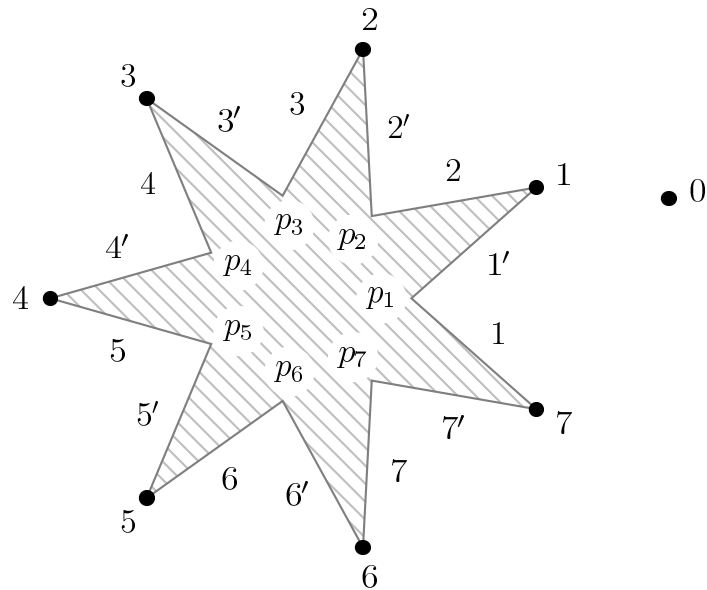}
\caption[Shape of the special sheet]{Shape of the special sheet (the
white region).} \label{fig:specjalneprzescieradlo}
\end{figure}

It follows that the special sheet looks as it is depicted on Figure
\ref{fig:specjalneprzescieradlo}, where points $p_1,p_2,\dots$ all
cover $\infty$ and those of them which carry the same labels should
be glued together. The pushing partition $\pi$ corresponding to the
pushing sequence $p_1,\dots,p_n$ can be therefore identified as a
receipt for gluing vertices $p_1,\dots,p_n$ of the special sheet
covering $\infty$. The meaning of the additional labeling of the
edges will be explained below.

\subsubsection{Collapsing non--special sheets}
\label{subsec:collapsing} Let us inflate the point $\infty$ a little
bit (or---speaking more precisely---while performing cuts from
points $1,2,\dots,n$ to $\infty$ let us leave some neighborhood of
$\infty$ untouched). The special sheet becomes now an oriented
surface with a boundary which coincides with a construction from
Section \ref{subsec:zaklejanie} where we glued a hole in a sphere
with the first collection of discs
(see also Figure 22 in \cite{Okounkov2000randompermutations}).

Let us collapse the union of non--special sheets (for details of
this construction we refer to
\cite{Okounkov2000randompermutations}). After this operation the
$2n$ edges of the special sheet (denoted by $1,1',2,2',\dots,n,n'$)
become glued together in pairs. We know that the vertices
$p_1,\dots,p_n$ are glued together according to the pushing
partition $\pi$;  this implies that the edges of the special sheet
are glued according to the pair partition $\pi_{\fat}$. This
provides an alternative description of the Okounkov's mapping $\Psi$
from coverings to maps on surfaces.

A careful reader will notice that the above statement is not
completely true: the reason is that Okounkov collapses non--special
sheets with two edges in a different way than the non--special
sheets with three or more edges. This means that in fact in the
process of computing $\pi_{\fat}$ we should treat two--element
blocks of $\pi$ differently. We leave the details as an exercise to
the reader.


\subsubsection{Many Jucys--Murphy elements}
In this article we deal with only one Jucys--Murphy element $J$
while in the work \cite{Okounkov2000randompermutations} one
considers products of many such elements $J_1,J_2,\dots$ and
therefore a careful reader might wonder if our methods are general
enough. However, one can easily translate Okounkov's statements
concerning products of many Jucys--Murphy elements into our favorite
language of products of moments of a single Jucys--Murphy element
which translates easily into properties of corresponding pushing
partitions.

\subsubsection{The case $(a_1 \xtra) \cdots (a_n \xtra)=e$}
Suppose now that $(a_1 \xtra)\cdots(a_n \xtra)=e$. It follows that
the construction from Section \ref{subsec:pierwszawzmiankaofppsq}
assigns to every cycle $b_s$ of the permutation $\pi_{\fat} \circ c$
the number $k_s=1$ (otherwise $(a_1 \xtra)\cdots(a_n \xtra)$ would
have a non--trivial cycle). Now we see from \eqref{eq:formulanakk}
that it is equivalent to the statement that every cycle $b_s$ of the
product $\pi_{\fat}\circ c$ makes exactly one counterclockwise wind
around the circle.

As in Section \ref{subsec:collapsing} let us collapse the
non--special sheets; the $2n$ vertices of the polygon representing
the special sheet fall into two classes: those which were covering
$\infty$ (and which correspond to glued vertices $p_1,\dots,p_n$)
and those corresponding to the ramification points $1,2,\dots,n$.
Suppose we are going clockwise around the boundary of the polygon
constituting the special sheet and thus we visit vertices in the
order $n,p_n,n-1,p_{n-1},\dots,1,p_1$; let us have a look in which
order we visit corners which meet in one of the vertices of the
second class. One can see that the edges which meet in such a vertex
correspond to one of the polygons which constitute the boundary of
the special sheet after gluing in $\infty$ or---in other words---to
one of the cycles of $\pi_{\fat}\circ c$. From our previous
discussion it follows that the corners will be visited in the
counterclockwise order; in Okounkov's notation this means that the
vertex is left.

Every group of the first class of vertices (covering $\infty$) which
is glued together corresponds to one of the blocks of the partition
$\pi$.
It follows that they will be visited in a clockwise order; in
Okounkov's notation this means that the vertex is right.

The above classification of vertices into the class of left and the
class of right ones is one of the key points of the paper
\cite{Okounkov2000randompermutations} and it is interesting that we
were able to reconstruct this result by the techniques of pushing
partitions.

\subsection[Connection with orthogonal polynomials]{Connection
with orthogonal polynomials and Wick product} Suppose that a
probability measure $\mu$ on the real line is given. We consider a
sequence $P_0,P_1,\dots$ of orthogonal polynomials with respect to
the measure $\mu$. Some questions concerning these polynomials are
extremely simple, for example the mean value is given by a trivial
equation $\int_\R P_k(x) d\mu(x)= \delta_{0k}$, but the price for
working with orthogonal polynomials is that we need a formula for
expressing a product $P_k P_l$ as a linear combination of
polynomials $P$ (in the non--commutative probability theory such
formulas are called Wick products). Also, a question arises how to
write monomials $x^k$ in this basis.

From this viewpoint one can easily see an analogy between orthogonal
polynomials and indicator functions $\Sigma$ considered in this
article. Some questions concerning $\Sigma$ are very simple: how to
evaluate its value on a given permutation or how to write a central
function as a linear combination of $\Sigma$. The price for this is
that we need to find formulas for computation of the products
$\Sigma_{k_1,\dots,k_m} \cdot \Sigma_{k'_1,\dots,k'_n}$ and need a
formula for expressing the moments of $J$ as a linear combination of
$\Sigma$. The `Wick product' \eqref{eq:definicjamnozenia} is very
similar to combinatorial objects which appear in this context in the
study of (generalized) Gaussian random variables \cite{GutaPhd}. It
would be very interesting to investigate the connections between
these two objects.

\section{Acknowledgments}
I thank Philippe Biane for introducing me into the subject and many
discussions. I also thank Roland Speicher, Ilona Kr\'olak and
Akihito Hora for many discussions and encouragement.

Research supported by State Committee for Scientific Research
(Komitet Bada\'n Naukowych) grant No.\ 2P03A00723; by EU Network
``QP-applications", contract HPRN-CT-2002-00729; by KBN-DAAD project
36/2003/2004. The research was conducted in \'Ecole Normale
Sup\'erieure du Paris, Institute des Hautes Etudes Scientifiques
(Bures--sur--Yvette, France), Syddansk Universitet (Odense, Denmark)
and Banach Center (Warszawa, Poland) on a grant funded by European
Post--Doctoral Institute for Mathematical Sciences.

%



\bibliographystyle{alpha}
\bibliography{biblio}


%


\end{document}